\documentstyle[11pt]{article}
\newtheorem{Lemma}{Lemma}
\newtheorem{Proposition}[Lemma]{Proposition}

\newtheorem{Theorem}[Lemma]{Theorem}

\newtheorem{Definition}[Lemma]{Definition}

\newtheorem{Corollary}[Lemma]{Corollary}

\newcommand{\lra}{\leftrightarrow}

\newcommand{\ed}{\ \stackrel{d}{=} \ }
\newcommand{\cd}{\ \stackrel{d}{\rightarrow} \ }
\newcommand{\cdlocal}{\ \rightarrow^d_{{\rm local}} \ }

\newcommand{\ef}{\stackrel{\rightarrow}{e}}

\newcommand{\eif}{\stackrel{\rightarrow}{e_1}}
\newcommand{\ejf}{\stackrel{\rightarrow}{e_j}}

\newcommand{\eiif}{\stackrel{\rightarrow}{e_2}}

\newcommand{\eb}{\stackrel{\leftarrow}{e}}

\newcommand{\WW}{\mbox{${\cal W}$}}
\newcommand{\GG}{\mbox{${\cal G}$}}
\newcommand{\FF}{\mbox{${\cal F}$}}

\newcommand{\MM}{\mbox{${\cal M}$}}

\newcommand{\eps}{\varepsilon}

\newcommand{\bT}{{\bf T}}
\newcommand{\bZ}{{\bf Z}}

\newcommand{\bV}{{\bf V}}
\newcommand{\bE}{{\bf E}}

\newcommand{\bc}{{\bf c}}

\newcommand{\bv}{{v}}
\newcommand{\by}{{y}}
\newcommand{\bw}{{\bf w}}

\newcommand{\var}{{\rm var}\ }

\newcommand{\sfrac}[2]{{\textstyle\frac{#1}{#2}}}
\newcommand{\PP}{{\cal P}}

\newcommand{\bW}{{\bf W}}

\newcommand{\dE}{\stackrel{\rightarrow}{\bE} }
\newcommand{\bm}{{\bf m}}
\newcommand{\bM}{{\bf M}}
\newcommand{\bmm}{\stackrel{\rightarrow}{\bm}}
\newcommand{\bMM}{\stackrel{\rightarrow}{\MM}}
\newcommand{\Mopt}{\MM_{{\rm opt}}}
\newcommand{\bMopt}{{\bMM}_{{\rm opt}}}

\newfam\Bbbfam
\font\tenBbb=msbm10
\font\sevenBbb=msbm7
\font\fiveBbb=msbm5
\textfont\Bbbfam=\tenBbb
\scriptfont\Bbbfam=\sevenBbb
\scriptscriptfont\Bbbfam=\fiveBbb

\title{The $\zeta(2)$ Limit in the Random Assignment Problem}

\author{David J. Aldous\thanks{Research supported by N.S.F. Grant DMS99-70901} \\
        Department of Statistics\\
        University of California\\
	367 Evans Hall \# 3860\\
        Berkeley CA 94720-3860\\
	aldous@stat.berkeley.edu\\
	http://www.stat.berkeley.edu/users/aldous}

\date{October 6, 2000}
\begin{document}
\maketitle

\begin{abstract}
The random assignment (or bipartite matching)
problem asks about 
$ A_n = \min_\pi \sum_{i=1}^n c(i,\pi(i))  $,
where $(c(i,j))$ is a $n \times n$ matrix with i.i.d. entries,
say with exponential(1) distribution, and the minimum is over
permutations $\pi$.
M\'{e}zard and Parisi (1987) 
used the replica method from statistical physics to argue
non-rigorously that $EA_n \to \zeta(2) = \pi^2/6$.
Aldous (1992) identified the limit in terms of a matching problem
on a limit infinite tree.
Here we construct the optimal matching on the infinite tree.
This yields a rigorous proof of the $\zeta(2)$ limit and of the
conjectured limit distribution of edge-costs and their rank-orders in the optimal matching.
It also yields the asymptotic essential uniqueness property:
every almost-optimal matching coincides with the optimal matching
except on a small proportion of edges.
\end{abstract}




\newpage
\section{Introduction}
\subsection{Background}
\begin{quote}
Consider the task of choosing an assignment of $n$ jobs to $n$
machines in order to minimize the total cost of performing
the $n$ jobs.
The basic input for the problem is an $n \times n$
matrix $(c(i,j))$, 
where $c(i,j)$ is viewed as the cost of performing job $i$
on machine $j$, and the {\em assignment problem} is to
determine a permutation $\pi$ that solves
\[ A_n = \min_\pi \sum_{i=1}^n c(i,\pi(i)) . \]
No doubt the simplest stochastic model for the assignment
problem is given by considering the $c(i,j)$ to be independent
random variables with the uniform distribution on $[0,1]$.
This model is apparently quite simple, but after some analysis
one finds that it possesses a considerable richness.
{\em Steele \cite{steele97}.}
\end{quote}
Part of this richness is that there have been four different approaches
to the study of the random assignment problem,
which we shall describe briefly.
See \cite{CS99,steele97} for further history.
Our focus is upon mathematical properties of $A_n$ and the
minimizing $\pi$, rather than algorithmic questions.

(a) {\em Rigorous bounds via linear programming}.
Steele (\cite{steele97} Chapter 4) goes on to give a detailed account of
results of Walkup \cite{wal79}, Dyer - Frieze - McDiarmid
\cite{DFM86} and Karp \cite{kar87} which lead to the upper bound
$ EA_n \leq 2  $.
The lower bound
$\limsup_n EA_n \geq 1 + e^{-1}$
was proved by Lazarus \cite{laz93} and 
subsequent work of
Olin \cite{olin} and
Goemans - Kodialam \cite{GK91} 
improved the lower bound to $1.51$.
Recently Coppersmith - Sorkin \cite{CS99} improved the
upper bound to $1.94$.

(b) {\em The replica method}.
M\'{e}zard - Parisi \cite{MP87} 
gave a non-rigorous argument based on the {\em replica method}
\cite{MPV87}
to show
$EA_n \to \pi^2/6 $.
The replica method and the somewhat related {\em cavity method} is extensively used in statistical physics
in non-rigorous analysis of spin glass and related disordered systems
\cite{parisi92}, and there is considerable interest in mathematical formalization
of the method
(see e.g. Talagrand
\cite{tal-ICM}).

(c) {\em Conjectured exact formulas}.
Parisi \cite{parisi98} conjectured that in the
case where the $c(i,j)$ have exponential($1$) distribution,
there is the exact formula
\begin{equation} EA_n = \sum_{i=1}^n i^{-2} . 
\label{exact}
\end{equation}
It is natural to seek an inductive proof.
Coppersmith - Sorkin \cite{CS99} and others \cite{AS99,BCR00,LW00} 
have recently formulated, and verified for small $n$,
more general conjectures concerning exact formulas for related problems, e.g. involving
incomplete matches on $m \times n$ bipartite graphs.
But (\ref{exact}) has not been proved rigorously
(an argument of Dotsenko \cite{dotsenko} is incomplete).

(d) {\em Weak convergence.}
Aldous \cite{me60} proved that
$\lim_n EA_n$ exists,
by identifying the limit as the value of a minimum-cost matching problem
on a certain random weighted infinite tree.

\subsection{Results}
In this paper we continue the ``weak convergence" analysis
mentioned above.   By studying the infinite-tree matching problem,
we obtain a rigorous proof of the M\'{e}zard - Parisi conclusion
$EA_n \to \pi^2/6$ (Theorem \ref{T1}), and also 
of their formula for the limit distribution of edge-costs in the 
optimal assignment (Theorem \ref{T2}),  
and of a conjecture of Houdayer et al \cite{HMM98}
concerning order-ranks of edges in the optimal assignment
(Theorem \ref{T4}).
Theorem \ref{T3} introduces and proves an
{\em asymptotic essential uniqueness} (AEU)
property of the random assignment problem:
roughly, AEU asserts that
every almost-optimal matching coincides with the optimal matching
except on a small proportion of edges.
Studying AEU is interesting for two reasons.
First, it can be defined for many ``optimization over random data"
problems (section \ref{sec-AEU}), providing a theoretical classification of such problems
(AEU either holds or fails in each problem)
somewhat in the spirit of computational complexity theory.
Second, in the statistical physics of disordered systems it
has been suggested \cite{MPV87} that the minima of the Hamiltonian 
should typically have an ultrametric structure, suggesting that in the
associated optimization problem the AEU property should fail.
Theorem \ref{T3} thus shows that the random assignment problem
(studied by physicists as a toy model of disordered systems)
has qualitatively different 
behavior than that predicted for more realistic models.

It follows easily from \cite{me60} that the distribution of the
$c(i,j)$ affects the limit of $EA_n$ only via the value of the
density function of $c(i,j)$ at $0$ (assuming this exists and is strictly
positive: see section \ref{sec-power} for the power-law case), and so in particular we may assume that the distribution
is exponential(1) instead of uniform.  
Indeed, to avoid later normalizations it is convenient to start out by
rephrasing the problem as
\[ A_n = \min_\pi \sfrac{1}{n} \sum_{i=1}^n c(i,\pi(i))  \]
where the $(c(i,j))$ are independent with exponential 
distribution with mean $n$. 
Write $\pi_n$ for the permutation attaining the minimum.

\begin{Theorem}
\label{T1}
$\lim_n EA_n = \pi^2/6$.
\end{Theorem}

\begin{Theorem}
\label{T2}
$c(1,\pi_n(1))$ converges in distribution; the limit distribution has
density
\[ h(x) =
\frac{e^{-x}(e^{-x}-1+x)}{(1-e^{-x})^2} ,
\ 0 \leq x < \infty . \]
\end{Theorem}

\begin{Theorem}
\label{T4}
For each $k \geq 1$ define
\[ q_n(k) = P(c(1,\pi_n(1))
\mbox{ is the $k$'th smallest of }
\{c(1,1), c(1,2),\ldots,c(1,n)\}). \]
Then $\lim_n q_n(k) = 2^{-k}$.
\end{Theorem}
\begin{Theorem}
\label{T3}
For each $0 < \delta < 1$
there exists $\eps(\delta) > 0$ such that,
if $\mu_n$ are permutations (depending on ($c(i,j)$)) such that
$E n^{-1} \#\{i: \mu_n(i) \neq \pi_n(i)\} \geq \delta$,
then
\[ \liminf_n E \left(
  \sfrac{1}{n} \sum_{i=1}^n c(i,\mu_n(i))  \right)
\geq \pi^2/6 + \eps(\delta) .\] 
\end{Theorem}
The logical structure of our proof is as follows.
\begin{itemize}
\item Review from \cite{me60} the infinite-tree minimum-weight matching problem 
and its connection with the finite-$n$ problem
(section \ref{sec-inf}).
\item Describe one matching $\Mopt$ in the infinite-tree setting (section \ref{sec-TC}).
\item Calculate the cost of the matching
$\Mopt$ and give an expression for the extra cost of any other matching
$\MM$ (section \ref{sec-PPr}).
\end{itemize}
In principle our method should yield an explicit lower bound for
$\eps(\delta)$ in Theorem \ref{T3} -- see 
section \ref{sec-quant}.

Matchings on the infinite tree are required to satisfy a certain
{\em spatial invariance} property, and this is the only
technically difficult ingredient of the paper.
Granted the abstract structure, the actual calculations
underlying the formulas in Theorems \ref{T1} - \ref{T4}
are quite straightforward, and we set out these calculations first in section
\ref{sec-calc}.
In a sense the calculation of the limit in Theorem \ref{T1}  
can be summarized in one sentence.
Let $X, X_i$ be i.i.d. with distribution determined by the identity
\begin{equation}
X \ed \min_i (\xi_i - X_i)
\mbox{ where $(\xi_i)$ is a Poisson($1$) process;}
\label{over1}
\end{equation}
then
\begin{equation}
\lim_n EA_n = \int_0^\infty x P(X_1+X_2 > x) \ dx .
\label{over2}
\end{equation}
From these it is straightforward (section \ref{sec-calc})
to evaluate the limit as $\pi^2/6$, though we do not have any calculation-free 
explanation of why such a simple number should appear
here or in Theorem \ref{T4}.

We remark that the technically  
most difficult, and conceptually most important,
part of the proofs of Theorems \ref{T1} - \ref{T3}
is the result from \cite{me60}
(restated as Theorem \ref{Told})
connecting the finite-$n$ assignment problem with the infinite-tree
problem.
The new additions of this paper are precisely
the construction (section \ref{sec-TC}) of the matching $\Mopt$
on the infinite tree, and its analysis in sections \ref{sec-analM}
and \ref{sec-PPr},
together with the elementary calculations in section \ref{sec-calc}.

\subsection{Discussion}
This paper was of course motivated by a desire to make the
M\'{e}zard - Parisi \cite{MP87} 
result rigorous.
In retrospect
our method, based on weak convergence to a limit random structure,
seems roughly similar to the non-rigorous
{\em cavity method} from statistical physics (see \cite{PM99,tal-ICM}
for brief descriptions), and Talagrand \cite{tal-SF} for its rigorous application in the Sherrington-Kirkpatrick model) though I do not understand that method well
enough to judge the degree of similarity.
So can our method be used for other mean-field models
of disordered systems in which the cavity method has been used?
It turns out it is easy to use our method to write down heuristic solutions
to three variations of the random assignment problem:
\begin{itemize}
\item The case where the density function $f_c(\cdot)$
of the costs $c(i,j)$ has $f_c(x) \sim x^r$ as $x \downarrow 0$
(section \ref{sec-power}).
\item The one-parameter Gibbs measure associated with the
assignment problem, recently studied rigorously by
Talagrand \cite{tal01} (section \ref{sec-Gibbs}).
\item The random combinatorial traveling salesman problem (section \ref{sec-TSP}).
\end{itemize}
We shall compare heuristic solutions from our method
with those from the cavity method.
Making a rigorous proof in these variations involves precisely the
issue in \cite{me60}: one needs a rigorous argument identifying
the optimal solution of an infinite-tree problem as the limit
of optimal solutions of finite-$n$ problems.

It is worth noting that in our method, calculations of limit
quantities depend on solving a problem-specific distributional
identity
(\ref{over1}, \ref{me2}, \ref{TSP1}, \ref{Gibbs-id}).
Such identities arise frequently in probabilistic analysis of
recursive algorithms \cite{RR92} and ``probability on trees"
\cite{lyonsPT96}, and seem worthy of further systematic study.

Section \ref{sec-FR} contains miscellaneous discussion;
the AEU property in more general contexts, 
the (lack of) insight our results cast on the exact conjecture
(\ref{exact}) and the intriguing analogy with Frieze's
$\zeta(3)$ result.

Our method rests upon a kind of ``local convergence" (\ref{eqlocal})
of $\pi_n$ to $\Mopt$, saying that the structure of $\pi_n$ and edge-costs
within a finite distance (interpreting edge-costs as distances) of a
typical vertex
converges in distribution to the structure of $\Mopt$ and edge-costs in the
infinite tree.
So in a sense it is Theorem \ref{T2} which is our basic result.
We can deduce Theorem \ref{T1} because
$EA_n = Ec(1,\pi_m(1))$, but we cannot study the variance of $A_n$ by our methods.
On the other hand we could write down expressions for the limit behavior of
arbitrarily complicated ``finite-distance" variations on Theorems \ref{T2}
and \ref{T4}, e.g. for
\begin{eqnarray*}
\lefteqn{\gamma(x,y,z):= \lim_n }\\ &&  \sfrac{1}{n}
E \left( \mbox{ number of pairs } (i,k): \ 
c(i,\pi(i)) \leq x, c(k,\pi(k)) \leq y, c(i,\pi(k)) \leq z
\right) . \end{eqnarray*}

\section{Calculations}
\label{sec-calc}
Consider the symmetric probability density
\[ f_X(x) = \left(e^{x/2} + e^{-x/2} \right)^{-2} , \quad - \infty < x < \infty
\]
which statisticians call the {\em logistic distribution}.
It is standard
\cite{logistic}
that the corresponding distribution function and variance are
\begin{equation}
F_X(x):= \int_{- \infty}^x f_X(y) dy
= \left( 1 + e^{-x}\right)^{-1}, \quad - \infty < x < \infty
\label{dF}
\end{equation}
\begin{equation}
\var X := \int_{-\infty}^\infty x^2 f_X(x) dx
= \pi^2/3
\label{var}
\end{equation}
and that the logistic distribution is characterized by the property
\begin{equation}
f(x) = F(x) (1-F(x)), \quad - \infty < x < \infty .
\label{char}
\end{equation}
In the calculations below we seek to exploit symmetry and
structure, rather than rely on ``brute force calculus".
\begin{Lemma}
\label{Llogistic}
Let
$0<\xi_1 < \xi_2 < \ldots$
be the points of a Poisson process of rate $1$.
Let $(X; X_i, i \geq 1)$ be independent random variables with
common distribution $\mu$.
Then
\begin{equation}
\min_{1 \leq i < \infty} (\xi_i - X_i) \ed  X
\label{min=d}
\end{equation}
if and only if $\mu$ is the logistic distribution.
\end{Lemma}
{\em Proof.}
The points $\{(\xi_i,X_i)\}$ form a Poisson point process
$\PP$ on $(0,\infty) \times (-\infty, \infty)$ with mean intensity
$\rho(z,x) dz dx = dz \ \mu(dx)$.
The distribution function $F$ of $\mu$ satisfies
\begin{eqnarray*}
1 - F(y) &=&
P\left(\min_{1 \leq i < \infty} (\xi_i - X_i) 
\geq y \right) \mbox{ by } (\ref{min=d})\\
&=& P\left(\mbox{no points of $\PP$ in $\{(z,x):\ z-x \leq y\}$}\right)\\
&=& \exp\left( - \int\int_{z-x\leq y} \ \rho(z,x) \ dzdx \right)\\
&=& \exp\left( - \int_0^\infty \bar{F}(z-y) \ dz \right), \mbox{ where }
\bar{F}(y) = 1-F(y) \\
&=& \exp\left( -\int_{-y}^{\infty} \ \bar{F}(u) du \right) .
\end{eqnarray*}
Differentiating, we see this is equivalent to
\begin{equation}
 F^\prime(y) = \bar{F}(-y)\bar{F}(y) . \label{Fsym}
\end{equation}
This implies the density $F^\prime(\cdot)$ is symmetric,
and so (\ref{Fsym}) is equivalent to the condition
(\ref{char})
characterizing the logistic distribution.
$\Box$

In some subsequent calculations we use the general identity
(for arbitrary real-valued r.v.'s $V,W$)
\begin{equation}
E(V-W)^+ = \int_{-\infty}^\infty P(V>x>W) \ dx .
\label{Fubini}
\end{equation}
Lemmas \ref{L1}, \ref{Lhform} and \ref{L2k} 
are the ``calculation" parts of Theorems \ref{T1}, \ref{T2}
 and \ref{T4} respectively.
\begin{Lemma}
\label{L1}
Let $X_1$ and $X_2$ be independent random variables with the logistic
distribution.  Then
\[ h(x) := P(X_1 + X_2 > x), \quad 0 \leq x < \infty \]
is the density of a probability distribution on $[0,\infty)$  with mean $\pi^2/6$.
\end{Lemma}
{\em Proof.}
To show it is a {\em probability} density,
\begin{eqnarray*}
\int_0^\infty h(x) dx &=& E(X_1+X_2)^+\\
&=&  E(X_1 - X_2)^+ \mbox{ by symmetry} \\
&=& \int_{-\infty}^\infty P(X_1 \geq y \geq X_2) \ dy
\mbox{ by } (\ref{Fubini})\\
&=& \int_{-\infty}^\infty
(1-F(y))F(y) \ dy\\
&=& \int_{-\infty}^\infty
f(y) \ dy \mbox{ by (\ref{char})}\\
&=& 1.
\end{eqnarray*}
And the mean is
\begin{eqnarray*}
\int_0^\infty x h(x) dx &=& \int_0^\infty x P(X_1 + X_2 \geq x) \ dx\\
&=& \sfrac{1}{2} E ((X_1+X_2)^+)^2\\
&=& \sfrac{1}{4} E(X_1+X_2)^2 \mbox{ by symmetry} \\
&=& \sfrac{1}{2} E X_1^2\\
&=& \pi^2/6
\mbox{ using (\ref{var})} .
\end{eqnarray*}
$\Box$

The next lemma derives an explicit formula for $h(x)$, though it
is not needed except to state the conclusion of Theorem \ref{T2}
(my thanks to Boris Pittel for this calculation).
\begin{Lemma}
\label{Lhform}
\begin{equation}  h(x) = 
\frac{e^{-x}(e^{-x}-1+x)}{(1-e^{-x})^2} ; 
\quad 0 \leq x < \infty .
\label{hformula}
\end{equation}
\end{Lemma}
{\em Proof.}
\begin{eqnarray*}
h(x) &=& \int_{- \infty}^\infty \frac{1}{(e^{u/2} + e^{-u/2})^2}
\ \frac{1}{1+e^{x-u}} \ du\\
&=& \int_{-\infty}^\infty \frac{e^u}{(e^u+1)^2} \ \frac{e^u}{e^x+e^u} \ du \\
&=& \int_0^\infty \frac{t}{(t+1)^2 (t+e^x)} \ dt\\
&=& \int_0^\infty \left[ \frac{e^x}{(e^x-1)^2}
\left(\frac{1}{t+1} - \frac{1}{t+e^x}\right)
- \frac{1}{e^x-1} \frac{1}{(t+1)^2} \right] \ dt \\
&=& \frac{e^x}{(e^x-1)^2} \left. \left( \log \frac{t+1}{t+e^x} \right)
  \right|_0^\infty 
\ + \ \frac{1}{e^x-1} \left. \frac{1}{t+1}
 \right|_0^\infty\\
&=& \frac{xe^x}{(e^x-1)^2} - \frac{1}{e^x-1}\\
\end{eqnarray*}
and this equals the formula at (\ref{hformula}).
$\Box$

\noindent
M\'{e}zard - Parisi (\cite{MP85} 
eq. (33))
write down the formula
\[ \tilde{h}(x) = \frac{x - e^{-x} \sinh x}{\sinh^2 x} \]
as the limit density for edge-costs in the optimal assignment,
in the essentially equivalent non-bipartite matching problem
on $2n$ vertices.
One can check
$\tilde{h}(x) = 2h(2x)$.
The factor of $2$ merely reflects the different normalization convention
(dividing by $2n$ instead of $n$).

\begin{Lemma}
\label{Lidentity2}
Let $(X,X_1,X_2,\eta)$ be independent, the $X$'s having logistic
distribution and $\eta$ having exponential($1$) distribution.
Then
\begin{equation}
X \ed \min(X_1,X_2) + \eta . \label{identity2}
\end{equation}
\end{Lemma}
{\em Proof.}
In the setting of Lemma \ref{Llogistic},
\[ (\xi_1,\xi_2,\xi_3,\ldots ) \ed (\eta, \eta+\xi_1^\prime, \eta+\xi_2^\prime, 
\ldots) \]
where the $(\xi_i^\prime)$ are a Poisson (rate $1$) process
independent of $\eta$.
So identity (\ref{min=d}) implies
\[  X \ed \eta + \min(-X_1,
\min_{i \geq 1}(\xi_i^\prime - X_{i+1})
) . \]
But by (\ref{min=d})
$\min_{i \geq 1}(\xi_i^\prime - X_{i+1})
\ed X_2$, so
\[ X \ed \eta + \min(-X_1,X_2) \]
and the result follows from the symmetry of the logistic distribution.
$\Box$

Note that (\ref{identity2}) does not characterize the logistic
distribution, because $X + $ constant will also satisfy (\ref{identity2}).
But it is not hard to show (Antar Bandyopadhyay, personal communication)
these are the only solutions.
\begin{Lemma}
\label{L2k}
In the setting of Lemma \ref{Llogistic}, for each $k \geq 1$
\[ \int_0^\infty P\left(x-X <
\min_{1 \leq i < \infty} (\xi_i - X_i) 
\mbox{ and }
x \geq \xi_k \right) \ dx = 2^{-k} . \]
\end{Lemma}
{\em Proof}.
Using symmetry of the logistic distribution, we may
rewrite the integrand as
$P(\xi_k \leq x < 
\min_{1 \leq i < \infty} (\xi_i + X_i) 
\ - X)$.
Then by (\ref{Fubini}) the value ($Q(k)$, say) of the integral is
\[ Q(k) = E \left(
\min_{1 \leq i < \infty} (\xi_i + X_i) 
- (\xi_k +X) \right)^+ . \]
Write
\[(\xi_1,\xi_2,\ldots) = (\xi_k - \eta_{k-1}, \xi_k - \eta_{k-2},
\ldots,\xi_k - \eta_1, \xi_k, \xi_k + \xi_1^\prime, \xi_k + \xi_2^\prime,
\ldots) , \]
so that $(\xi^\prime_i)$ is a Poisson process.
Then
\begin{eqnarray*}
Q(k) \! &=&
E\left(
\min(X_1 - \eta_{k-1},\ldots,X_{k-1} - \eta_1,X_k, 
\min_{1 \leq i < \infty} (\xi^\prime_i + X_{k+i}) 
) \ -X \right)^+ \\
&=& E(
\min(X_1 - \eta_{k-1},\ldots,X_{k-1} - \eta_1,X_k, 
X_{k+1}
) \ -X)^+ \\
&& \mbox{ by (\ref{min=d}) and symmetry }\\
&=& \! \int_{-\infty}^\infty \! P(X < x <
\min(X_1 - \eta_{k-1},\ldots,X_{k-1} - \eta_1,X_k, X_{k+1})
) dx \mbox{ by } (\ref{Fubini})\\
&=& \int_{-\infty}^\infty 
F(x) P(x<
\min(X_1 - \eta_{k-1},\ldots,X_{k-1} - \eta_1,X_k))
\ (1-F(x)) dx\\
&=& \int_{-\infty}^\infty 
f(x) P(x<
\min(X_1 - \eta_{k-1},\ldots,X_{k-1} - \eta_1,X_k))
\ dx \mbox{ by } (\ref{char})\\
&=& P(X<
\min(X_1 - \eta_{k-1},\ldots,X_{k-1} - \eta_1,X_k)) .
\end{eqnarray*}
Since $(\eta_1,\ldots\eta_{k-1})$ are distributed as the first $k-1$ points
of a Poisson process, we can relabel variables to obtain
\[ Q(k) = P(X< \min(X_0,X_1-\xi_1,X_2-\xi_2,\ldots,X_{k-1}-\xi_{k-1})) . \]
For $k=1$ this says
$Q(1) = P(X<X_0) = 1/2$,
so it is enough to show $Q(k) = \frac{1}{2}Q(k-1)$ for $k \geq 2$.
Write $M = \min(X,X_0)$.
Then $P(X=M) = 1/2$ and this event is independent of the value of $M$,
so
\[ Q(k) = \sfrac{1}{2} P(M< \min(X_1-\xi_1,X_2-\xi_2,\ldots,X_{k-1}-\xi_{k-1})) . \]
Writing
\[ (\xi_1,\xi_2,\xi_3,\ldots ) \ed (\eta, \eta+\xi_1^\prime, \eta+\xi_2^\prime, 
\ldots) \]
and setting $X^\prime_i = X_{i+1}$ gives
\[ Q(k) = \sfrac{1}{2} P(M + \eta < \min(X^\prime_0,X^\prime_1-\xi^\prime_1,\ldots,
X^\prime_{k-2} - \xi^\prime_{k-2})) . \]
But Lemma \ref{Lidentity2} shows $M+\eta \ed X$ and so we have
proved $Q(k) = \sfrac{1}{2} Q(k-1)$.

\newpage
\section{Working on the infinite tree}
\subsection{The infinite matching problem}
\label{sec-inf}
In this section we review results from Aldous \cite{me60}, which
unfortunately require considerable space to describe.

Let $\bV$ be the set of finite words $\bv = v_1v_2\ldots v_d, \ 0 \leq d <
\infty$
where each $v_i$ is a natural number; include in $\bV$ the
empty word $\phi$.
There is a natural tree $\bT$ with vertex-set $\bV$ and edge-set $\bE$, where an
edge $e \in \bE$ is of the form $e = (\bv, \bv j), j \geq 1$,
where for
$\bv = v_1v_2 \ldots v_d$
we write $\bv j = v_1v_2 \ldots v_d j$ for the $j$'th {\em child}
of $\bv$, and call $\bv$ the {\em parent} of $\bv j$.
Now attach random edge-weights as follows.
For each $\bv \in \bV$ let the weights 
$(W(\bv,\bv j), j \geq 1)$ on the
edges $((\bv, \bv j), j \geq 1)$
be the points 
of a Poisson (rate $1$) point process on $(0,\infty)$,
independently as $\bv$ varies.
Call this structure the {\em Poisson-weighted infinite tree} (PWIT).
See figure 1 for an illustration.
Write $\lambda$ for the probability distribution of the whole configuration
$(W(e))$ of edge-weights.
So $\lambda$ is a probability measure on the space
$\bW = (0,\infty)^{\bE}$
of all possible configurations $\bw = (w(e), e \in \bE)$ of edge-weights.

\setlength{\unitlength}{0.5in}
\begin{picture}(9.5,4.8)(-0.5,0)
\put(0,0){\fbox{11}}
\put(1,0){\fbox{12}}
\put(2,0){\fbox{13}}
\put(3,0){\fbox{21}}
\put(4,0){\fbox{22}}
\put(5,0){\fbox{23}}
\put(6,0){\fbox{31}}
\put(7,0){\fbox{32}}
\put(8,0){\fbox{33}}
\put(1,2){\fbox{1}}
\put(4,2){\fbox{2}}
\put(7,2){\fbox{3}}
\put(4,4){\fbox{$\phi$}}
\put(1.15,1.75){\line(0,-1){1.33}}
\put(0.95,1.75){\line(-1,-2){0.68}}
\put(1.35,1.75){\line(1,-2){0.68}}
\put(4.15,1.75){\line(0,-1){1.33}}
\put(3.95,1.75){\line(-1,-2){0.68}}
\put(4.35,1.75){\line(1,-2){0.68}}
\put(7.15,1.75){\line(0,-1){1.33}}
\put(6.95,1.75){\line(-1,-2){0.68}}
\put(7.35,1.75){\line(1,-2){0.68}}
\put(4.15,3.75){\line(0,-1){1.33}}
\put(3.9,3.83){\line(-3,-2){2.47}}
\put(4.4,3.83){\line(3,-2){2.47}}
\put(1.9,2.9){0.5}
\put(3.92,2.9){0.8}
\put(6.05,2.9){2.1}
\put(0.1,0.9){1.8}
\put(0.92,0.9){1.9}
\put(1.88,0.9){3.7}
\put(3.1,0.9){0.9}
\put(3.92,0.9){1.2}
\put(4.87,0.9){2.5}
\put(6.1,0.9){0.4}
\put(6.92,0.9){2.8}
\put(7.87,0.9){4.1}
\end{picture}

\vspace{0.1in}

{\bf Figure 1.}
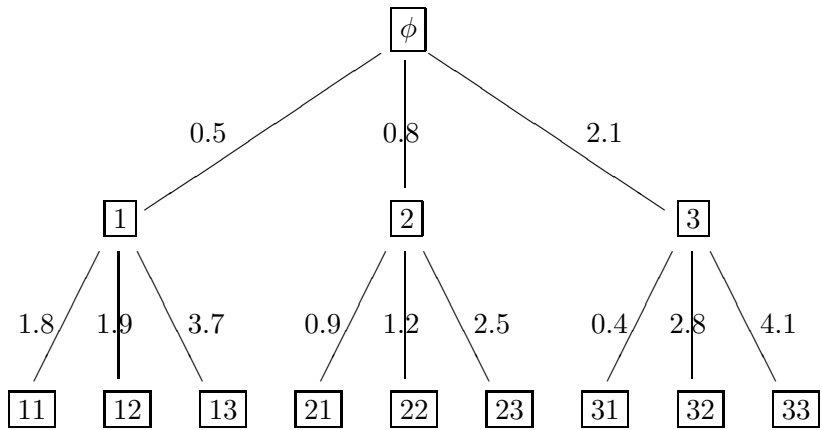
Part of a realization $\bw$ of the PWIT.  
The weight $w(e)$ is written next to the edge $e$.

\vspace{0.2in}
\noindent
A {\em matching} $\bm$ on the PWIT is a set of edges of $\bT$ such that each vertex
is incident to exactly one edge in the set.
Formally we can identify the set $\bM$ of matchings as the subset
$\bM \subseteq \{0,1\}^{\bE}$ defined by:
$\bm = (m(e))\in \bM$ iff $\sum_{e: \bv \in e} m(e) = 1 \ \forall \bv \in \bV$.
A matching can also be regarded as a map from vertices to vertices,
in which case we write $\bmm (\bv) = \bv^\prime$ to indicate $(\bv,\bv^\prime)$ is
an edge in the matching.
See figure 2 for an illustration.

\setlength{\unitlength}{0.5in}
\begin{picture}(9.6,6.8)(-0.5,-2)
\put(0,0){\fbox{11}}
\put(1,0){\fbox{12}}
\put(2,0){\fbox{13}}
\put(3,0){\fbox{21}}
\put(4,0){\fbox{22}}
\put(5,0){\fbox{23}}
\put(6,0){\fbox{31}}
\put(7,0){\fbox{32}}
\put(8,0){\fbox{33}}
\put(1,2){\fbox{1}}
\put(4,2){\fbox{2}}
\put(7,2){\fbox{3}}
\put(4,4){\fbox{$\phi$}}
\put(0.95,1.75){\line(-1,-2){0.68}}
\put(6.95,1.75){\line(-1,-2){0.68}}
\put(4.15,3.75){\line(0,-1){1.33}}
\put(3.9,2.9){0.8}
\put(0.1,0.9){1.8}
\put(6.1,0.9){0.4}
\put(2.6,-2){\fbox{213}}
\put(3.05,-0.25){\line(-1,-6){0.22}}
\put(3.9,-2){\fbox{221}}
\put(4.15,-0.25){\line(0,-1){1.33}}
\put(5.35,-0.25){\line(1,-6){0.22}}
\put(5.4,-2){\fbox{232}}
\put(2.46,-1.1){2.7}
\put(3.9,-1.1){1.1}
\put(5.54,-1.1){2.1}
\end{picture}

\vspace{0.1in}

{\bf Figure 2.}
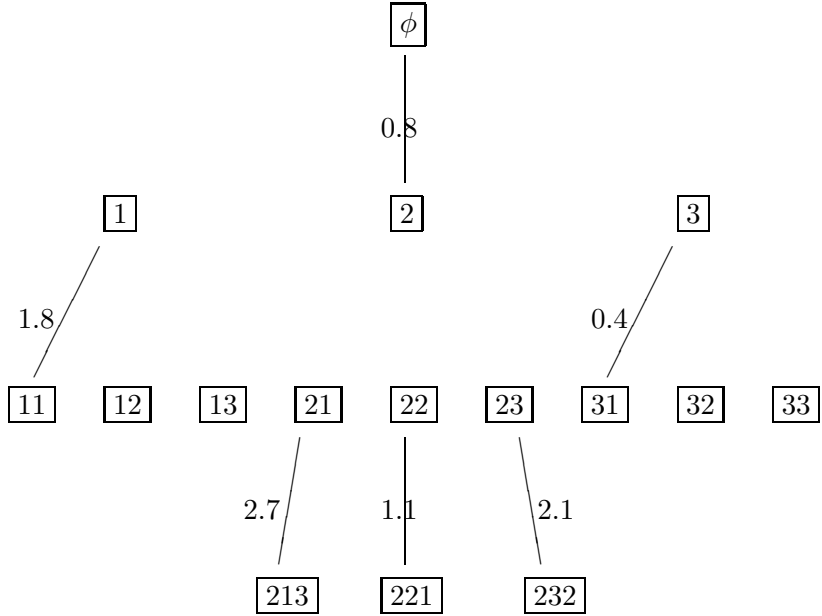
Part of a matching $\bm$ on the PWIT.  
Only the edges of $\bm$ are drawn.

\vspace{0.2in}
\noindent
We study random matchings $\MM$ on the PWIT.
Such a random matching will be dependent on the edge-weights,
so formally we are dealing with a probability measure $\mu$ on
$\bW \times \bM$ (the joint distribution of edge-weights and indicators
of edges in the matching) with marginal distribution $\lambda$ on $\bW$.

The connection between matchings $\MM$ on the PWIT and matchings $\pi_n$
in the $n \times n$ random assignment problem is via
{\em local convergence} of matchings on $\bT$ induced by the
{\em unfolding map}; these ideas are explained at 
the end of this section.
Theorem \ref{Told} below says (roughly speaking)
that $\lim_nEA_n$ is the average cost per edge in a minimum-cost
matching on $\bT$.
At first sight such a result looks false, because there is
a greedy matching $\MM_{{\rm greedy}}$ consisting of edges
$(\phi, 1), \ (2,21), \ (3,31), \ldots, (11,111), \ 
(12,121), \ldots \ldots$,
and $\MM_{{\rm greedy}}$ has average cost equal to $1$,
which is certainly not the desired $\lim_n EA_n$.
But the precise result is more complicated, because 
matchings in Theorem \ref{Told} are required to satisfy a certain property:
informally
\begin{quote}
(*) The rule for whether an edge $e$ is in the matching should be spatially
invariant, that is should not depend on which vertex was chosen as the root
of $\bT$.
\end{quote}
We now work towards making this notion precise.
Our definitions are superficially different from those in \cite{me60},
but we reconcile them in section \ref{sec-rec}.

For each $\bw \in \bW$ and each $i \geq 1$ we define an isomorphism of $\bT$
(that is, a bijection $\bV \to \bV$ which induces a bijection $\bE \to \bE$;
we denote either bijection by $\theta_i^\bw(\cdot)$)
which expresses the idea
``make vertex $i$ the root, and relabel vertices to preserve order structure".
Precisely, first define $k \in \{1,2,3,\ldots \}$ by
\[ w(i,i(k-1)) < w(\phi,i) < w(i,ik) \]
interpreting the left as $0$ for $k=1$.
Then define

$\theta_i^\bw(i) = \phi$

$\theta_i^\bw(\phi) = k$

$\theta_i^\bw$ takes vertices $(i1,i2,\ldots,i(k-1); ik, i(k+1),\ldots)$ to
vertices $(1,2,\ldots, k-1; k+1, k+2, \ldots)$

for $\bv = i j_2 \ldots j_l$ and $l \geq 3$ let
$\theta_i^\bw (\bv) = (\theta_i^\bw(i j_2)) j_3 \ldots j_l$

for $\bv = j_1j_2 \ldots j_l$ where $j_1 \neq i$ let
$\theta_i^\bw(\bv) = kj_1j_2 \ldots j_l$.

\vspace{0.07in}
\noindent
As already mentioned, the space $\bW \times \bM$ describes possible
combinations $(\bw,\bm)$ of edge-weights $\bw$ and a matching $\bm$.
Consider the enlarged state space
$\bZ = \bW \times \bM \times \{1,2,3,\ldots \}$ with elements $(\bw,\bm,k)$
where we interpret the $k$ as meaning that vertex $k$ is distinguished.
The maps $\theta_i^\bw$ induce a single map
$\theta: \bZ \to \bZ$ which we interpret as
``relabel the distinguished vertex as the root, relabel other vertices to
preserve order structure, and make the previous root into the distinguished
vertex".  Precisely,
$\theta(\bw,\bm,i) = (\hat{\bw},\hat{\bm},k)$ where

$w(e) = \hat{w}(\theta_i^\bw(e))$

$m(e) = \hat{m}(\theta_i^\bw(e))$

$k = \theta_i^\bw(\phi)$.\\
See figure 3.
Note $\theta = \theta^{-1}$.

A probability measure $\mu$ on $\bW \times \bM$
(describing a random matching on the PWIT)
extends to a $\sigma$-finite measure
$\mu^* := \mu \times {\tt count}$ on $\bZ$,
where {\tt count} is counting measure on $\{1,2,3,\ldots \}$.
\begin{Definition}
\label{def-si}
A random matching $\MM$ on the PWIT, with distribution $\mu$ (say) on
$\bW \times \bM$, is {\rm spatially invariant} if 
$\mu^*$ is invariant under $\theta$, that is if
$\mu^*(\cdot) = \mu(\theta^{-1}(\cdot))$
as $\sigma$-finite measures on $\bZ$.
\end{Definition}

\setlength{\unitlength}{0.5in}
\begin{picture}(9.5,6.8)(-0.5,-2)
\put(0,0){\fbox{11}}
\put(1,0){\fbox{12}}
\put(2,0){\fbox{13}}
\put(3,0){\fbox{21}}
\put(4,0){\fbox{22}}
\put(5,0){\fbox{23}}
\put(6,0){\fbox{31}}
\put(7,0){\fbox{32}}
\put(8,0){\fbox{33}}
\put(1,2){\fbox{1}}
\put(4,2){\fbox{2}}
\put(7,2){\fbox{3}}
\put(4,4){\fbox{$\phi$}}
\put(-0.4,-2){\fbox{111}}
\put(0.9,-2){\fbox{121}}
\put(2.4,-2){\fbox{132}}
\put(0.05,-0.25){\line(-1,-6){0.22}}
\put(1.15,-0.25){\line(0,-1){1.33}}
\put(2.35,-0.25){\line(1,-6){0.22}}
\put(1.15,1.75){\line(0,-1){1.33}}
\put(0.95,1.75){\line(-1,-2){0.68}}
\put(1.35,1.75){\line(1,-2){0.68}}
\put(4.15,1.75){\line(0,-1){1.33}}
\put(3.95,1.75){\line(-1,-2){0.68}}
\put(4.35,1.75){\line(1,-2){0.68}}
\put(7.15,1.75){\line(0,-1){1.33}}
\put(6.95,1.75){\line(-1,-2){0.68}}
\put(7.35,1.75){\line(1,-2){0.68}}
\put(4.15,3.75){\line(0,-1){1.33}}
\put(3.9,3.83){\line(-3,-2){2.47}}
\put(4.4,3.83){\line(3,-2){2.47}}
\put(1.9,2.9){0.8}
\put(2.4,2.6){*}
\put(3.9,2.9){0.9}
\put(6.05,2.9){1.2}
\put(0.1,0.9){0.5}
\put(0.9,0.9){2.1}
\put(1.9,0.9){3.2}
\put(3.1,0.9){2.2}
\put(3.9,0.9){2.6}
\put(4.9,0.9){2.7}
\put(4.75,0.5){*}
\put(6.1,0.9){1.1}
\put(6.5,0.57){*}
\put(6.9,0.9){1.7}
\put(7.9,0.9){3.3}
\put(-0.54,-1.1){1.8}
\put(0,-1.05){*}
\put(0.9,-1.1){0.4}
\put(0.93,-0.91){*}
\put(2.54,-1.1){1.6}
\put(2.3,-1.05){*}
\end{picture}

\vspace{0.1in}

{\bf Figure 3.}
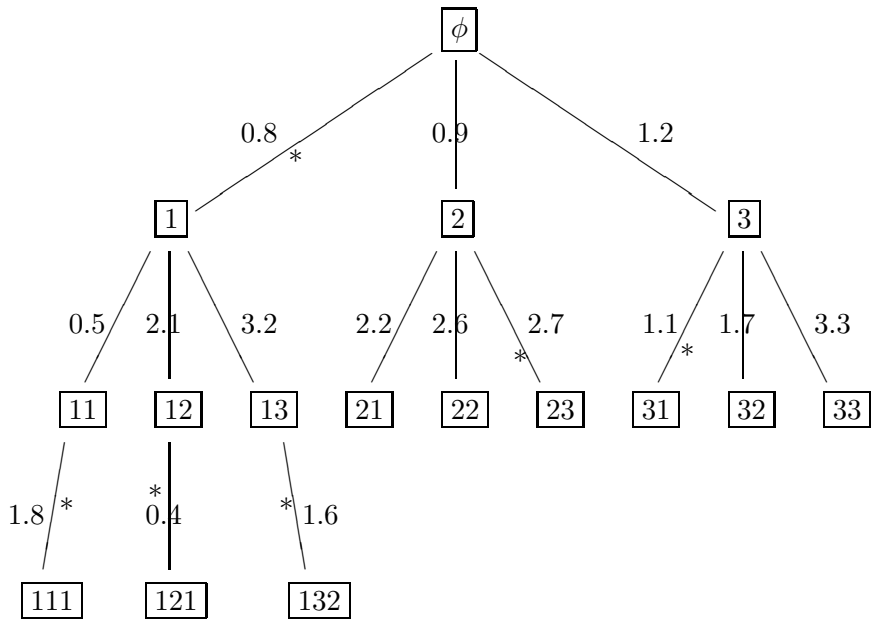
Take the edge-weights $\bw$ in figure 1, and the matching $\bm$ in figure 2,
and distinguish vertex $2$.  Then figure 3 illustrates part of $\theta(\bw,\bm,2) 
= (\hat{\bw},\hat{\bm},1)$.  Figure 3 shows the new edge-weights $\hat{w}$ attached to the
edges; edges in the new matching $\hat{\bm}$ are indicated by $*$; the
new distinguished vertex is vertex $1$.

\vspace{0.2in}
\noindent
In a random matching $\MM$ the root $\phi$ is matched to some
random neighbor vertex $\bMM(\phi) \in \{1,2,3,\ldots\}$,
and the edge $(\phi, \bMM(\phi))$ has some random cost
$W(\phi,\bMM(\phi))$.
We can now reformulate the main result of \cite{me60} as
\begin{Theorem}
\label{Told}
\[
\lim_n EA_n =   \inf 
\{ E W(\phi, \bMM(\phi)) \} \] 
where the infimum is over all spatially invariant random matchings $\MM$
on the PWIT.
\end{Theorem}
Though Theorem \ref{Told} is the only result from \cite{me60} needed to prove
Theorem \ref{T1}, the other new theorems need further infrastructure from
\cite{me60}.  The reader should perhaps skip the rest of this section on first
reading.

A random matching in the $n \times n$ random assignment problem
(in brief, an $n$-matching)
was denoted by $\pi_n$ in the introduction, but can be reformulated
as a $n \times n \ \{0,1\}$-valued random matrix $\MM = (m(i,j))$.
Call an $n$-matching {\em spatially invariant} if
the joint distribution $((c(i,j),m(i,j)), 1 \leq i,j \leq n)$
is invariant under the automorphisms of the complete bipartite graph,
that is under permutations of $i$, under permutations of $j$,
and under complete interchange of $i$ and $j$.
Given any random $n$-matching, by applying a uniform random
automorphism one obtains a spatially invariant $n$-matching with
the same distribution of average-cost-per-edge $A_n$, so there
is no loss of generality in considering only spatially invariant
$n$-matchings.
There is a notion of ``local convergence"
$\MM_n \cdlocal \MM$, made precise at (\ref{eqlocal}) later,
which implies in particular that
$c(1,\bMM_n(1)) \cd W(\phi,\bMM(\phi))$.
The following two results (stated slightly differently in \cite{me60} -- see section \ref{sec-rec}) immediately imply
Theorem \ref{Told}.
\begin{Theorem}
\label{Told-hard}
Let $\MM$ be a spatially invariant matching on the PWIT with
$E W(\phi, \bMM(\phi))  
< \infty$.
Then there exist spatially invariant $n$-matchings $\MM_n$ such that
$\MM_n \cdlocal \MM$
and
$(c(1,\bMM_n(1)), \ 1 \leq n < \infty)$ is uniformly intergrable.
\end{Theorem}
\begin{Theorem}
\label{Told-easy}
Let $j_n \uparrow \infty$, and let $\MM_{j_n}$ be spatially invariant
$j_n$-matchings such that
$\limsup_n 
E c(1,\bMM_{j_n}(1)) < \infty$. 
Then there exists a subsequence $k_n$ of $j_n$ such that
$\MM_{k_n} \cdlocal \MM$, where $\MM$ is some spatially invariant
matching on the PWIT.
\end{Theorem}
{\em Remarks.}
The proof of Theorem \ref{Told-easy} is fairly simple, using
compactness arguments.
The proof of Theorem \ref{Told-hard} in \cite{me60} is difficult
and lengthy.
In light of the new results of this paper, we only need Theorem
\ref{Told-hard} for $\MM = 
\Mopt$, and perhaps the explicit structure of $\Mopt$
could be used to simplify the proof of Theorem \ref{Told-hard}.

To formalize {\em local convergence} requires further notation.
Fix $n$ and let $\bV^{(n)}$ be the set of vertices
$\bv = v_1v_2 \ldots v_l$
with $v_1 \leq n$ and $v_i \leq n-1$ for $i \geq 2$.
Write $\bT^{(n)} = (\bV^{(n)},\bE^{(n)})$ for the corresponding
subtree of $\bT$.
Let $G_{nn}$ be 
the complete bipartite graph, with vertex-set 
$\{1,2,\ldots,n\} \times \{1,2\}$.
Given a realization $\bc = (c(i,j))$ of the cost matrix,
one can define a graph homomorphism $\psi = \psi_{\bc}$ from $\bT^{(n)}$ onto $G_{nn}$
as follows.

$\psi(\phi) = (1,1)$

for $i \in \{1,2,3,\ldots, n\},
$ define $
\psi(i) = (j,2)$, for the $j$ such that
$c(1,j)$ is the $i$'th smallest of 
$\{c(1,u), 1 \leq u \leq n\}$

for $i \in \{1,2,3,\ldots, n-1\}$
and $i^\prime \in \{1,2,\ldots,n\}$, define $
\psi(i^\prime i) = (k,1)$, for the $k$ such that
$c(k,1)$ is the $i$'th smallest of 
$\{c(u,\psi(i^\prime)), 1 \leq u \leq n,  u \neq 1\}$

and then inductively: for $\bv = v_1v_2\ldots v_{2m}$,
for $i \in \{1,2,3,\ldots, n-1\},
$ define $
\psi(\bv i) = (j,2)$, for the $j$ such that
$c(\psi(\bv),j)$ is the $i$'th smallest of 
$\{c(\psi(\bv),u), 1 \leq u \leq n, \ (u,2) \neq\psi(\bv^-)\}$
where $\bv^-$ is the parent of $\bv$

for $\bv = v_1v_2 \ldots v_{2m+1}$,
for $i \in \{1,2,3,\ldots, n-1\},
$ define $
\psi(\bv i) = (k,1)$, for the $k$ such that
$c(k,\psi(\bv))$ is the $i$'th smallest of 
$\{c(u,\psi(\bv)), 1 \leq u \leq n, \ (u,1) \neq\psi(\bv^-)\}$.

This {\em folding} map $\psi_{\bc}$ induces an {\em unfolding}
map which uses the matrix $\bc$ to define edge-weights $(W^{(n)}(e))$ on the
edge-set $\bE^{(n)}$:
\begin{eqnarray*}
W^{(n)}(\bv, \bv k) = c(i,j) && \mbox{ if }
\psi_{\bc}(\bv) = (i,1) \mbox{ and } \psi_{\bc}(\bv k) = (j,2)\\
&& \mbox{ or if }
\psi_{\bc}(\bv) = (j,2) \mbox{ and } \psi_{\bc}(\bv k) = (i,1).
\end{eqnarray*}
Now for $h \geq 1$ define $\bE_{(h)} \subset \bE$ as the set of edges
each of whose end-vertices is of the form
$v_1v_2 \ldots v_l$ with $l \leq h$ and $\max_i v_i \leq h$.
So $\bE_{(h)}$ is a finite set of edges.
It is straightforward to see
that for fixed $h$
\[ (W^{(n)}(e), e \in \bE_{(h)}) \cd
(W(e), e \in \bE_{(h)}) \]
where the limits are the edge-weights in the PWIT.
(Essentially, this is the fact that the order statistics
of $n$ independent exponential (mean $n$) r.v.'s converge
to the points of a Poisson (rate $1$) process.)
Recall that we represent a random $n$-matching as $\{0,1\}$-valued
random variables $\MM_n(e)$ indexed by the edges $e$ of $G_{nn}$.
Use the homomorphism $\psi$
(considered as a map on edges)
to define
\[ \widetilde{\MM}_n(e) = \MM_n(\psi(e)), \ e \in \bE^{(n)} . \]
We can now define
{\em local convergence}
$\MM_n \cdlocal \MM$ to mean:
for each fixed $h$,
\begin{equation}
\left( (W^{(n)}(e),\widetilde{\MM}_n(e)), \ e \in \bE_{(h)} \right)
\cd
\left((W(e),\MM(e)), \ e \in \bE_{(h)} \right) . 
\label{eqlocal}
\end{equation}

\subsection{Heuristics for the construction of $\Mopt$}
\label{sec-heur}
Underlying the rigorous construction in the next section is a simple
heuristic idea.
Given a realization of the PWIT, consider defining
\begin{eqnarray}
 \quad \quad \quad X_{\phi} &=&
\mbox{ cost of optimal matching on } \bT\nonumber\\
 &-& \mbox{ cost of optimal matching on } \bT \setminus \{\phi\} . \label{eq-heur}
\end{eqnarray}
Here we mean {\em total} cost, so we get
$\infty - \infty$, and so (\ref{eq-heur}) makes no sense as a
rigorous definition, though statistical physics uses such
renormalization arguments all the time.
But pretend (\ref{eq-heur}) does make sense.
Then for each $\bv \in \bV$ we can 
define $X_{\bv}$ similarly in terms of the subtree $\bT^{\bv}$
rooted at $\bv$ (the vertices of $\bT^{\bv}$ are $\bv$ and its descendants):
\begin{eqnarray*}
 \quad \quad \quad X_{\bv} &=&
\mbox{ cost of optimal matching on } \bT^{\bv}\\
 &-& \mbox{ cost of optimal matching on } \bT^{\bv} \setminus \{\phi\} .
\end{eqnarray*}
And we get the recursion
\begin{equation}
 X_{\bv} = \min_{1 \leq j < \infty} 
(W(\bv,\bv j) - X_{\bv j}) \label{rec-heur}
\end{equation}
because the left side is the cost difference between using or
not using $\bv$ in a matching on $\bT^{\bv}$; to use an edge $(\bv,\bv j)$
we have to pay the cost of the edge and the difference between
the cost of not using or using $\bv j$, which is the right side.
Moreover in the optimal matching on $\bT^{\bv}$, vertex $\bv$ should
be matched to the vertex $\bv j$ attaining the minimum in (\ref{rec-heur}).

In the next section we show that one can make a rigorous argument
by first constructing (by fiat) random variables satisfying the 
recursion (\ref{rec-heur}), then using these random variables to define
a matching.  Not having interpretation (\ref{eq-heur}) 
means it's not obvious rigorously that 
this matching is optimal, but it turns out (see start of proof of Proposition
\ref{P9}) that {\em weak} optimality is quite easy to prove.

\subsection{The construction}
\label{sec-TC}
Each edge $e \in \bE$ of $\bT$ corresponds to two directed edges
$\ef, \eb$: write $\dE$ for the set of directed edges.
For directed edges we have the language of family relationships: each
edge $\ef = (\bv^\prime,\bv)$ has an infinite number of children 
of the form $(\bv,\by), \ \by \neq \bv^\prime)$.
Thus the directed edge $(273,27)$ has children
$(27,2), \ (27,271), \ (27,272), \ (27,274), \ldots $.
The PWIT has edge-weights $W(e)$ on undirected edges, and we write
$W(\ef) = W(\eb) = W(e)$ for the directed edges $\ef, \eb$ corresponding to $e$.
\begin{Lemma}
\label{Ljoint}
Jointly with the edge-weights $(W(e), e \in \bE)$
of the PWIT we can construct
$\{X(\ef), \ef \in \dE\}$
such that \\
(i) each $X(\ef)$ has the logistic distribution\\
(ii) for each $\ef$, with children $\eif,\eiif,\ldots$ say, 
\begin{equation}
X(\ef) = \min_{1 \leq j < \infty} (W(\ejf) - X(\ejf)) \label{Xrec}
\end{equation}
\end{Lemma}
{\em Proof.}
For a vertex $\bv = i_1i_2 \ldots i_h$ write $|\bv| = h$.
For $h \geq 1$ write
\begin{eqnarray*}
\dE_h &=& \{\ef = (\bv, \bv j): \ |\bv| = h, j \geq 1 \}\\
\dE_{\leq h} &=& \{\ef = (\bv,\by): \ |\bv| \leq h, |\by| \leq h \}.
\end{eqnarray*}
Take independent logistic random variables $\{X(\ef): \ef \in \dE_h\}$,
independent of the family $(W(e))$.
Then use (\ref{Xrec}) recursively to define $X(\ef)$ for each 
$\ef \in \dE_{\leq h}$.
Lemma \ref{Llogistic} and the natural independence structure ensures that
each $X(\ef)$ has logistic distribution.
The construction specifies a joint distribution for
$\{W(e), e \in \bE; \ X(\ef), \ef \in \dE_h \cup \dE_{\leq h}\}$.
The Kolmogorov consistency theorem completes the proof.
$\Box$

{\em Remarks.}
(a) By modifying on a null set, we shall assume that the minimum
in (\ref{Xrec}) is attained by a unique $j$.\\
(b) For an edge $\ef = (v^\prime,v)$ directed away from $\phi$,
the $X(\ef)$ constructed above formalizes the notion of
$X_v$ in the previous section.\\
(c) 
The same type of construction can be associated with general
fixed-point identities for distributions which have an appropriate
format.  See e.g. Aldous \cite{me89} for its use in studying
a model of ``frozen percolation" on the infinite binary tree.

The corollary below (whose routine proof we omit) formalizes the independence structure
implicit in the construction.
For a directed edge $\ef$, consider the set consisting of
all its descendant edges, but not $\ef$ itself; 
then write $D(\ef)$ for this set of edges,
considered as {\em undirected} edges.
\begin{Corollary}
\label{Cindep}
For each $i \geq 1$ let
$f_{i,1},f_{i,2},\ldots \in \bE$ and let 
$\stackrel{\to}{e_{i,1}},
\stackrel{\to}{e_{i,2}},
\ldots \in \dE$.
Write $D_i = \cup_j D(
\stackrel{\to}{e_{i,j}})
\ \cup \{f_{i,1},f_{i,2},\ldots \}$.
If the sets $D_i$ are disjoint
as $i$ varies
then the $\sigma$-fields
$\sigma(X(
\stackrel{\to}{e_{i,j}}), W(f_{i,j}), j \geq 1)$
are independent
as $i$ varies.
\end{Corollary}
{\em Remark.}
It is natural to ask whether $X(\ef)$ is $\sigma (W(e), e \in D(\ef))$-measurable,
in other words whether the influence of the ``boundary" $\dE_h$ in the Lemma
\ref{Ljoint} construction vanishes in the $h \to \infty$ limit.
In the terminology of statistical physics, this asks about extremality
of free boundary Gibbs states.  We have not studied this question carefully,
but there are reasons to suspect the answer is negative.

For $v, v^\prime \in \bV$ write $v \sim v^\prime$ if
$(v,v^\prime)$ is an undirected edge.
Fix a realization of 
$(W(e), e \in \bE; \ X(\ef), \ef \in \dE)$.
For each $v \in \bV$ define
\begin{equation}
v^* = \arg \min_{v^\prime \sim v} (W(v, v^\prime) - X(v,v^\prime)) 
\label{defM0}
\end{equation}
that is,
$v^*$ is the vertex attaining the minimum.
\begin{Lemma}
\label{LcM}
The set of undirected edges $\{(v,v^*): \ v \in \bV\}$
is a matching on the PWIT.
\end{Lemma}
{\em Proof.}
It suffices to check that $(v^*)^* = v$ for each $v \in \bV$.
Fix $v$.  Then
\begin{eqnarray*}
W(v,v^*) - X(v,v^*) &<& \min_{y \sim v, y \neq v^*} 
(W(v,y) - X(v,y))
\mbox{ by definition of } v^*\\
&=&X(v^*,v) \mbox{ by } (\ref{Xrec}) .
\end{eqnarray*}
In other words
\begin{equation}
X(v,v^*) + X(v^*,v) > W(v,v^*) . \label{XXW}
\end{equation}
Suppose $(v^*)^* = z \neq v$; then
\[ W(v^*,z) - X(v^*,z) < W(v^*,v) - X(v^*,v) . \]
But by (\ref{Xrec})
\[ X(v,v^*) \leq W(v^*,z) - X(v^*,z) . \]
Adding these last two inequalities gives
$X(v,v^*) < W(v^*,v) - X(v^*,v)$
which contradicts (\ref{XXW}) and establishes the lemma.
$\Box$

Write $\Mopt$ for the random matching specified by
Lemma \ref{LcM}.
Note that the argument leading to (\ref{XXW}) can be reversed,
to give a more symmetric criterion for whether an edge is in $\Mopt$:
\begin{equation}
e \mbox{ is an edge of $\Mopt$ iff }
W(e) < X(\ef) + X(\eb)
\label{Mchar}
\end{equation}
where $\ef, \eb$ are the directed edges corresponding to $e$.
It seems intuitively clear that $\Mopt$ should be spatially invariant;
we verify this in section \ref{sec-SI} as Lemma \ref{LMSI}.

{\em Remarks.}
(d) It is not obvious whether $\Mopt$ is a function of the
weights $(W(e), e \in \bE)$ only, or involves additional external randomization.
In fact by (\ref{Mchar}) this question is equivalent to the question in remark (c) above.

(e) Intuitively, the quantities $X(\ef)$ play a role analogous to 
{\em dual variables} in linear programming.  But
we are unable to make this idea more precise.

\subsection
{Analysis of $\Mopt$}
\label{sec-analM}
Consider the random cost
$ W(\phi, \bMopt(\phi)) $ 
of the edge $(\phi,\bMopt (\phi))$
containing the root $\phi$ in $\Mopt$.
\begin{Proposition}
\label{P8}
(a) The random variable
$ W(\phi, \bMopt(\phi)) $ 
has the probability density function $h(\cdot)$ described in Lemma \ref{L1},
and so
$  EW(\phi, \bMopt(\phi)) 
= \pi^2/6$.\\
(b) 
\[ P(\bMopt(\phi) = k) = 2^{-k}, \ k \geq 1 . \]
\end{Proposition}
{\em Proof.}
The edge-weights $(W(\phi,i), i \geq 1)$
and the $X$-values $(X(\phi,i), i \geq 1)$
are distributed as the Poisson process $(\xi_i)$
and the i.i.d. logistics $(X_i)$ in Lemma \ref{Llogistic}.
Using the latter notation and the definition (\ref{defM0}) of
$\Mopt$,
\[
 W(\phi, \bMopt(\phi))  
= \xi_I, \mbox{ where }
I = \arg \min_{i \geq 1} (\xi_i - X_i) . \]
As a sophisticated way to obtain this distribution
(there are alternate, elementary ways)
fix $0<y<\infty$ and condition on the event
$A_y := \{\exists J: \xi_J = y\}$.
Conditionally, the other points $(\xi_i, i \neq J)$ and other
$X$-values $(X_i, i \neq j)$ are distributed as a Poisson
process $(\xi^\prime_j, j \geq 1)$ say and i.i.d. logistics
$(X^\prime_j, j \geq 1)$ say, independent of $X_J$,
whose conditional distribution remains logistic.
So
\[ P(I=J |A_y) = P(y - X_J <X^\prime)
\mbox{ where }
X^\prime = \min_{j \geq 1} (\xi^\prime_j - X^\prime_j) . \]
But by Lemma \ref{Llogistic} $X^\prime$ has logistic distribution;
since it is independent of $X_J$, we see from the definition of $h(\cdot)$
that
$P(I=J|A_y) = h(y)$.
Then
\[ P(\xi_I \in [y,y+dy]) = P(I=J|A_y) 
P(\mbox{ some $\xi_i$ in 
$[y,y+dy]$}) = h(y) dy \]
establishing part (a).  For part (b) we need to show
\begin{equation}
P(I \geq k+1) = 2^{-k}, \ k \geq 1 .
\label{Ik}
\end{equation}
From the argument above
\[ P(I=J, I \geq k+1|A_y)
= P(y-X_J< \min_{j \geq 1} (\xi^\prime_j - X^\prime_j), \ 
\xi^\prime_k < y)
= h_k(y), \mbox{ say } . \]
So
\[ P(\xi_I \in [y, y+dy], I \geq k+1)
= h_k(y) dy \]
and then
\[ P(I \geq k+1) = \int_0^\infty h_k(y) dy \]
and the integral is evaluated by Lemma \ref{L2k}.

\begin{Proposition}
\label{P9}
Let $\MM$ be a spatially invariant matching of the PWIT such that
$P(\bMM(\phi) \neq \bMopt(\phi)) > 0$.
Then
$  EW(\phi, \bMM(\phi)) 
>  EW(\phi, \bMopt(\phi)) $.
\end{Proposition}
Before embarking upon the proof of Proposition \ref{P9}, let us show how
Theorems \ref{T1} - \ref{T3} are deduced from 
Propositions \ref{P8} and \ref{P9}.
Indeed, Propositions \ref{P8} and \ref{P9}
imply
\[ \inf
\{ E W(\phi, \bMM(\phi))   
: \ \bMM \mbox{ is spatially invariant matching on the PWIT} \}
= \sfrac{\pi^2}{6} \]
and so Theorem \ref{Told} implies Theorem \ref{T1}.
We next claim that the 
optimal assignments in the $n$-matching problem,
considered as random matchings $\MM_n$ as in section \ref{sec-inf},
satisfy
\begin{equation}
\mbox{
$\MM_n \cdlocal \Mopt$
and
$(c(1,\bMM_n(1)), \ 1 \leq n < \infty)$ is uniformly intergrable.
} \label{cmc}
\end{equation}
If not, then by Theorem \ref{Told-easy} some subsequence
converges locally to some $\MM \neq \Mopt$, so 
$\limsup_n EA_n \geq 
 E W(\phi, \bMM(\phi))   
> \pi^2/6$,
the last inequality by Proposition \ref{P9},
but this contradicts Theorem \ref{T1}.
Now (\ref{cmc}) implies Theorems \ref{T2} and \ref{T4},
because the definition of local convergence (\ref{eqlocal})
implies
\[ \left(\bMM_n(1),c(1,1),c(1,2),\ldots \right)
\cd \left( \bMopt(\phi), W(\phi,1), W(\phi,2), \ldots \right) \]
and then Proposition \ref{P8}
identifies the required limit distributions.
To prove Theorem \ref{T3}, fix $\delta > 0$.
Suppose the assertion were false for that $\delta$;
then we can find random
$j_n$-matchings $\MM^\prime_{j_n}$
(for some $j_n \to \infty$)
such that
\begin{equation}
 P\left( \bMM^\prime_{j_n}(1) 
\neq \bMM_{j_n}(1) \right) \geq \delta
\label{za} \end{equation}  where as above $\MM_{j_n}$ is the optimal $j_n$-matching; and
\begin{equation}
E c\left(1,\bMM^\prime_{j_n}(1) \right) \to \pi^2/6 . \label{zb}
\end{equation}
But by Theorem \ref{Told-easy} we can, after passing to a further
subsequence, assume $\MM^\prime_{j_n} \cdlocal \MM$, 
for some spatially invariant random matching $\MM$.
Then (\ref{za}) implies 
$P(\bMM(\phi) \neq \bMopt(\phi)) \geq \delta$,
while (\ref{zb}) implies 
$  EW(\phi, \bMM(\phi)) 
\leq  EW(\phi, \bMopt(\phi)) $,
contradicting Proposition \ref{P9}.
This ``proof by contradiction" establishes 
Theorem \ref{T3}.

Examining the argument above, we see that a weak inequality
in Proposition \ref{P9} would be sufficient to prove Theorem \ref{T1},
while the strict inequality is needed for our proofs of
Theorems \ref{T2} - \ref{T3}.

\subsection{Proof of Proposition 18.}
\label{sec-PPr}
Let $\MM$ be a spatially invariant matching, and write
$A$ for the event 
$\{\bMM(\phi) \neq \bMopt(\phi)\}$.
On $A$ write 
\[ (v_{-1},v_0,v_1) = (\bMopt(\phi),\phi,\bMM(\phi)). \]
Then on $A$ we can define a doubly-infinite path
$\ldots, v_{-2},v_{-1},v_0,v_1,v_2,\ldots$
by: $\forall - \infty < m < \infty$
\[ (v_{2m-1},v_{2m}) \mbox{ is an edge of } \Mopt \]
\[ (v_{2m},v_{2m+1}) \mbox{ is an edge of } \MM .\]
In section \ref{sec-SI} we shall use spatial invariance to prove
\begin{Lemma}
\label{L7}
Conditional on $A$, 
the distributions of
$X(v_{-2},v_{-1})$
and
$X(v_{0},v_{1})$
are the same.
\end{Lemma}
Now 
\begin{eqnarray}
X(v_{-2},v_{-1})&=&
\min_{y \sim v_{-1},y \neq v_{-2}}
(W(v_{-1},y) - X(v_{-1},y))
\mbox{ by (\ref{Xrec})}\nonumber\\
&=&
(W(v_{-1},v_0) - X(v_{-1},v_0))
\label{12}
\end{eqnarray}
 by (\ref{defM0}), because
$(v_{-1},v_0)$ is in $\Mopt$.
Also, by (\ref{Xrec})
\begin{equation}
X(v_{-1},v_0) \leq
W(v_{0},v_1) - X(v_{0},v_1)
\label{01}
\end{equation}
so that
\[ D:= W(v_0,v_1) - X(v_0,v_1) - X(v_{-1},v_0) \geq 0 . \]
Combining (\ref{12}) with this definition of $D$ gives
\[ W(v_0,v_1) - W(v_0,v_{-1})
= D + X(v_0,v_1) - X(v_{-2},v_{-1}) . \]
So
\begin{eqnarray}
\lefteqn{EW(\phi, \bMM(\phi)) 
-  EW(\phi, \bMopt(\phi)) 
}\nonumber\\
&=& E(
 W(v_0,v_1) - W(v_0,v_{-1})
)1_A \nonumber\\
&=&
ED1_A 
+ EX(v_0,v_{-1})1_A - EX(v_{-2},v_{-1})1_A \nonumber\\
&=& ED1_A
\mbox{ by Lemma \ref{L7}}. \label{DA}
\end{eqnarray}
Since $D \geq 0$ this is enough to establish the weak
inequality
$  EW(\phi, \bMM(\phi)) 
\geq  EW(\phi, \bMopt(\phi)) $.
To establish the strict inequality, suppose 
that to the contrary
$  EW(\phi, \bMM(\phi)) 
=  EW(\phi, \bMopt(\phi)) $.
Then (\ref{DA}) implies $ED1_A = 0$, and then 
(\ref{01}) implies that on $A$
we have
$X(v_{-1},\phi) =
W(\phi,v_1) - X(\phi,v_1)$.
And this implies
that on $A$ we have
\[
 v_1 = \arg \min_i \, \! ^{[2]} (W(\phi,i) - X(\phi,i))
\]
where $\min^{[2]}$ denotes the second-smallest value,
because $\bMM(\phi) = v_{-1} = \arg \min_i (W(\phi,i) - X(\phi,i))$.
So without restricting to $A$ we have
\[
P \! \left( \!  \bMM(\phi) = 
    \arg \min_i (W(\phi,i) - X(\phi,i))
\mbox{ or }  \arg \min_i \, \! ^{[2]} (W(\phi,i) - X(\phi,i))
\! \right) = 1 .
\]
By an immediate use of spatial invariance we must have the same
property at every $v \in \bV$:
\begin{equation}
P \left( \bMM(v) = 
    \arg \min_i (W(v,i) - X(v,i))
\mbox{ or }  \arg \min_i \, \! ^{[2]} (W(v,i) - X(v,i))
\right) = 1 .
\label{M12}
\end{equation}
So it's enough to show this can't happen.
\begin{Proposition}
\label{Lonly}
The only spatially invariant random matching on the PWIT
satisfying (\ref{M12}) is $\Mopt$.
\end{Proposition}
{\em Remark.}
One can regard Proposition \ref{Lonly} as a kind of
``subcritical percolation" fact.
There is some set
${\MM}_{{\rm 2-opt}}$
of edges $(v^\prime, v^{\prime \prime})$ such that
\[
 v^{\prime \prime} = \arg \min_{v \sim v^\prime} \, \! ^{[2]} (W(v^\prime,v) - X(v^\prime,v))
\mbox{ and }
 v^{\prime} = \arg \min_{v \sim v^{\prime \prime}} \, \! ^{[2]} (W(v^{\prime \prime},v) - X(v^{\prime \prime},v))
 . \]
In contrast to $\Mopt$, the first equality here does not imply the second.
One can show that Proposition \ref{Lonly} is equivalent to the assertion that
there is no infinite path consisting of alternating edges from
$\Mopt$ and
${\MM}_{{\rm 2-opt}}$.

{\em Proof.}
Define a path $\phi = w_0,w_{-1},w_{-2},w_{-3},\ldots$
inductively by: for $m = 1,2, \ldots$
\begin{eqnarray*}
w_{-2m+1} &=& \arg \min_{y \sim w_{-2m+2}} (W(w_{-2m+2},y)-X(w_{-2m+2},y))\\
w_{-2m} &=& \arg \min_{y \sim w_{-2m+1}} \, \! ^{[2]}(W(w_{-2m+1},y)-X(w_{-2m+1},y)) .
\end{eqnarray*}
Then
\[ (w_{-2m+1},w_{-2m+2}) \mbox{ is an edge of } \Mopt , 
\ \mbox{each } m = 1,2,3,\ldots . \]
Suppose $\MM$
were a spatially invariant matching satisfying (\ref{M12})
such that $A:= 
\{\bMM(\phi) \neq \bMopt(\phi)\}$
has $P(A) > 0$.
Then on $A$ we have
\[ (w_{-2m},w_{-2m+1}) \mbox{ is an edge of } \MM , 
\ \mbox{each } m = 1,2,3,\ldots . \]
(The construction of $(w_m)$ resembles the previous construction of
$(v_m)$, but note that the definitions of $(w_m)$ above and
of $B_m, \bar{B}_q, B^*$ below 
involve only the
PWIT and the $(X(\ef))$, not any hypothetical matching $\MM$.)
So by (\ref{M12}) for $v = w_{-2m}$
\[ A \subseteq B_m :=
\{w_{-2m+1} = \arg \min_{y \sim w_{-2m}} \, \! ^{[2]}(W(w_{-2m},y)-X(w_{-2m},y)) \} 
\]
and so \[A \subseteq B:= \cap_{m=1}^\infty B_m . \]
Writing $\bar{B}_q := \cap_{m=1}^q B_m$ we have
\begin{equation}
 \mbox{ if } P(A) > 0 \mbox{ then }
\lim_{q \to \infty} P(\bar{B}_{q+1})/P(\bar{B}_q) = 1 . \label{PABQ}
\end{equation}
In section \ref{sec-SI} we shall use spatial invariance to
prove
\begin{Lemma}
\label{LBB}
$P(\bar{B}_{q+1}) = P(\bar{B}_q \cap B^*)$, where
\[ B^* := \{
 \phi = \arg \min_{y \sim w_1} \, \! ^{[2]} (W(w_1,y) - X(w_1,y))
\} \]
where
\[ w_1 := \arg \min_{y \sim \phi} \, \! ^{[2]} (W(\phi,y) - X(\phi,y))
 . \]
\end{Lemma}
Using Lemma \ref{LBB},
\[ \lim_{q \to \infty} \frac{P(\bar{B}_{q+1})}{P(\bar{B}_q) } = 
\lim_q \frac{P(\bar{B}_q \cap B^*)}{P(\bar{B}_q)}
= \frac{P(B \cap B^*)}{P(B)}
= P(B^*|B) . \]
So by (\ref{PABQ}), to prove Proposition \ref{Lonly} and hence
Proposition \ref{P9} it is enough to prove
\begin{Lemma}
\label{LBB*}
If $P(B)>0$ then $P(B^*|B)<1$. 
\end{Lemma}
{\em Remark.}
Intuitively, Lemma \ref{LBB*} seems clear for the following reason.
Event $B$ depends only on what happens on the branch from $\phi$
through $w_{-1}$, while $B^*$ depends only on what happens
on the branch through $w_1$.
Even though there is dependence between the branches, the dependence
shouldn't be strong enough to make $B^*$ conditionally 
{\em certain} to happen.
Formalizing this idea requires study of the effect of conditioning the
PWIT, and we defer completing the proof
until section \ref{sec-CI2},
though the next lemma is one ingredient of the proof.

\begin{Lemma}
\label{LXI}
Define
\begin{eqnarray*}
X^\downarrow &=& \min_{i \geq 1} (W(\phi,i) - X(\phi,i))\\
I &=& \arg \min_{i \geq 1} (W(\phi,i) - X(\phi,i)).
\end{eqnarray*}
For $-\infty < b < a < \infty$ define
\[ g(a,b) = P \left(
W(\phi,I) -b >
\left. 
  \min_{k \geq 1} \, \! ^{[2]}  
(W(I,Ik) - X(I,Ik)) \ \right| X^\downarrow = a \right) . \]
Then $g(a,b) > 0$.
\end{Lemma}
{\em Proof.}
Note that (\ref{Xrec}) shows
\[ X(\phi,I) = \min_{k \geq 1} (W(I,Ik)-X(I,Ik)) . \]
Conditionally on $X(\phi,I) = x$,
the other values of
$\{W(I,Ik)-X(I,Ik), \ k \geq 1\}$
are the points of a certain inhomogeneous Poisson process
on $(x,\infty)$, and so
\begin{equation}
P \left(
\left. 
  \min_{k \geq 1} \, \! ^{[2]}  
(W(I,Ik) - X(I,Ik)) \in [y,y+dy] \ \right| X(\phi,I) =x  \right) 
= \beta_x(y) dy
\label{beta1}
\end{equation}
for a certain function $\beta_x(\cdot)$ such that
$\beta_x(y) > 0$ for all $y>x$.
The quantities in (\ref{beta1}) are independent of $W(\phi,I)$
and so (\ref{beta1}) remains true if we also condition on
$W(\phi,I) = a+x$.
So
\begin{eqnarray*}
\lefteqn{ \tilde{g}(a,b,x) :=}\\
&&P \left(
\left. 
W(\phi,I) -b >
  \min_{k \geq 1} \, \! ^{[2]}  
(W(I,Ik) - X(I,Ik)) \ \right| 
X(\phi,I) = x, \ W(\phi,I) = a+x \right) 
\end{eqnarray*}
satisfies
$\tilde{g}(a,b,x) > 0$ for all $-\infty < b < a < \infty$
and $- \infty < x < \infty$.
Since $X^\downarrow = W(\phi,I) - X(\phi,I)$ we see
\[ g(a,b) = E(\tilde{g}(a,b,X(\phi,I)) | X^\downarrow = a) > 0 \]
as required.

\section{Some technical details}
\subsection{Spatial invariance proofs}
\label{sec-SI}
The definition (Definition \ref{def-si}) of spatial invariance for a 
random matching $\MM$ on the PWIT
involves the (probability) distribution $\mu$ of
$((\MM(e),W(e)), \ e \in \bE)$
and the 
$\sigma$-finite measure
$\mu^* := \mu \times {\tt count}$ on $\bZ$.
Intuitively, if we use an arbitrary rule to distinguish
some vertex $k$ (the rule depending on the realization of
$\MM$ and the $W(e)$),
then the $\mu^*$-distribution of the resulting configuration
on $\bZ$ is just the $\mu$-distribution of the configuration without
the vertex being distinguished.

To see the point of spatial invariance, let us show that the
greedy matching
$\MM_{{\rm greedy}}$
from section \ref{sec-inf} is not spatially invariant.

\setlength{\unitlength}{0.5in}
\begin{picture}(9.5,5.9)(-0.5,-0.3)
\put(0,0){\fbox{11}}
\put(0,2){\fbox{1}}
\put(0,4){\fbox{$\phi$}}
\put(1.5,0){\fbox{12}}
\put(1.5,2){\fbox{2}}
\put(0.15,1.8){\line(0,-1){1.43}}
\put(0.15,3.8){\line(0,-1){1.43}}
\put(0.44,1.8){\line(3,-4){1.08}}
\put(0.44,3.8){\line(3,-4){1.08}}
\put(-0.3,0.8){0.8}
\put(-0.3,2.8){1.5}
\put(0.7,0.8){2.1}
\put(0.7,2.8){2.4}
\put(5,0){\fbox{1}}
\put(5,2){\fbox{$\phi$}}
\put(5,4){\fbox{2}}
\put(6.5,0){\fbox{3}}
\put(6.5,2){\fbox{21}}
\put(5.15,1.8){\line(0,-1){1.43}}
\put(5.15,3.8){\line(0,-1){1.43}}
\put(5.44,1.8){\line(3,-4){1.08}}
\put(5.44,3.8){\line(3,-4){1.08}}
\put(4.7,0.8){0.8}
\put(4.7,2.8){1.5}
\put(5.7,0.8){2.1}
\put(5.7,2.8){2.4}
\put(0.7,4.7){$\bw \in B^*$}
\put(5.7,4.7){$\theta^{
-1}(\bw) \in \tilde{B}$}
\end{picture}

\vspace{0.1in}

{\bf Figure 4.}
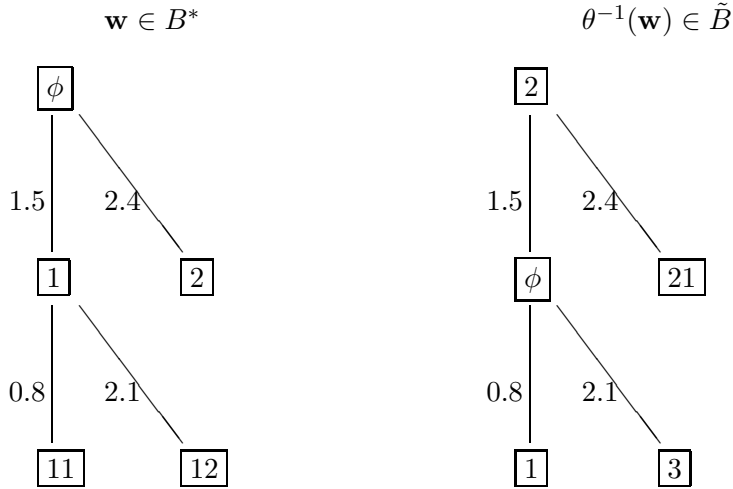
Part of a realization of $\bw \in B^*$ and of $\theta^{-1}(\bw) \in \tilde{B}$.

\vspace{0.13in}
\noindent
Consider the event
$B:= \{W(1,11)<W(\phi,1)<W(1,12)\}$
which has probability $1/4$.
So
$B^*:= B \cap \{ 1 \mbox{ is distinguished } \}$
has $\mu^*$-measure $1/4$.
The inverse image $\tilde{B} = \theta^{-1}(B^*)$ is
(see figure 4)
\[ \tilde{B} = \{W(\phi,2) < W(2,21), \ 2 \mbox{ is distinguished} \} \]
and this has $\mu^*$-measure $1/4$ because
$P(W(\phi,2)<W(2,21)) = 1/4$.
(This equality of $\mu^*$-measures is a consequence of the fact that the
distribution of edge-weights is spatially invariant.)
Now for a random matching $\MM$ to be spatially invariant, the
event
\[ B^* \cap \{\bMM(\phi) = 1\} \]
must have the same $\mu^*$-measure as its inverse image under $\theta$,
which is the event
\[ \tilde{B} \cap \{\bMM(\phi) = 2 \} . \]
But for the greedy matching we always have $\bMM(\phi) = 1$,
so the former event has $\mu^*$-measure $1/4$
while the latter event has $\mu^*$-measure $0$.

{\em Proof of Lemma \ref{L7}.}
In words, the idea is that the distribution of
$X(v_0,v_1)$ as seen from $\phi = v_0$ is the same
(by spatial invariance) as the distribution of
$X(v_0,v_1)$ as seen from $v_2$, 
which (just by relabeling) is the distribution of
$X(v_{-2},v_{-1})$ as seen from $v_0 = \phi$. 
Precisely:
\[
P(A, X(v_0,v_1) \in \cdot) =
P\left(\bMM(\phi) \neq \bMopt(\phi),
X(\phi,\bMM(\phi)) \in \cdot \right)\]
\[ = \mu^* \left(\bMM(\phi) \neq \bMopt(\phi),
X(\phi,\bMM(\phi)) \in \cdot, 
\bMM(\phi) \mbox{ is distinguished} \right)
\]
\[ = \mu^* \left(\bMM(\phi) \neq \bMopt(\phi),
X(\bMM(\phi),\phi) \in \cdot,
\bMM(\phi) \mbox{ is distinguished} \right)
\]
using spatial invariance to switch the root from $\phi$ to $\bMM(\phi)$,
and noting that the event
$\{\bMM(\phi) \neq \bMopt(\phi)\}$
is identical to the event
$\{\bMM(\bMM(\phi)) \neq \bMopt(\bMM(\phi))\}$;
\[ = \mu^* \left(\bMM(\phi) \neq \bMopt(\phi),
X(\bMM(\phi),\phi) \in \cdot, 
\bMopt(\phi) \mbox{ is distinguished} \right)
\]
because $\mu^*$-measure doesn't depend on the rule for distinguishing
a vertex;
\[ = \mu^* \left(\bMM(\bMopt(\phi)) \neq \bMopt(\bMopt(\phi)),
X(\bMM(\bMopt(\phi)),\bMopt(\phi)) \in \cdot, 
\bMopt(\phi) \mbox{ is distinguished} \right)
\]
using spatial invariance to switch the root from $\phi$ to $\bMopt(\phi)$;
\[ = \mu^* \left(\bMM(\phi) \neq \bMopt(\phi),
X(\bMM(\bMopt(\phi)),\bMopt(\phi)) \in \cdot, 
\bMopt(\phi) \mbox{ is distinguished} \right)
\]
because the events are identical;
\begin{eqnarray*}
 &=& P \left(\bMM(\phi) \neq \bMopt(\phi),
X(\bMM(\bMopt(\phi)),\bMopt(\phi)) \in \cdot 
\right) \\
&=& P(A, X(v_{-2},v_{-1}) \in \cdot )
\end{eqnarray*}
because $v_{-1} = \bMopt(\phi)$ and
$v_{-2} = \bMM(v_{-1})$.

{\em Proof of Lemma \ref{LBB}.}
This is essentially the same argument as above.
In words, the idea is that the probability of
$\bar{B}_q \cap B^*$ as seen from $\phi = w_0$ is the same
(by spatial invariance) as the probability of
$\bar{B}_q \cap B^*$ as seen from $w_2 
= \arg \min_{w \sim w_1} (W(w_1,w)-X(w_1,w))$,
which (just by relabeling) is the probability of
$\bar{B}_{q+1}$ as seen from $w_0 = \phi$.

We leave details to the reader.

\begin{Lemma}
\label{LMSI}
$\Mopt$ is spatially invariant.
\end{Lemma}
{\em Proof.}
$\Mopt$ is determined by the $W(e)$ and the $X(\ef)$,
and the $X(\ef)$ satisfy the deterministic relation (\ref{Xrec})
which is unaffected by relabeling vertices.
Moreover the joint distribution of the $X(\ef)$ is determined by
the fact (Lemma \ref{Ljoint}) that the r.v.'s
$(X(\ef), \ef \in \dE_h)$
are independent logistic.
In words, we need to show that the property
``the $(X(\ef), \ef \in \dE_h)$
are independent logistic"
is preserved under $\theta$.

Fix $k$ and write
\[ B = \{(X(\ef), \ef \in \dE_h) \in C, \ k \mbox{ is distinguished}\}  \]
for arbitrary $C$ in the appropriate range space.
Write
\[ A_l = \{W(l,k-1) < W(\phi,l) < W(l,k) \} . \]
Then (see figure 5)
\[ \theta^{-1}(B ) =
\cup_l \left[ A_l \cap \{ l \mbox{ is distinguished }\} \cap
\{(X(\ef), \ef \in \dE_{h,l}) \in C\} \right] \]
where $\dE_{h,l}$ is the set of edges $(v,vj)$ with $j \geq 1$
and with $v = j_1j_2 \ldots j_{h-1}$ for $j_1 \neq l$,
or with $v = lj_1j_2 \ldots j_h$.
Since $\sum_l P(A_l) = 1$ and the $X(\ef)$ under consideration
are independent of $A_l$, verifying
$\mu^*(B) = \mu^*(\theta^{-1}(B))$
reduces to showing that
$\{X(\ef), \ef \in \dE_{h,l} \}$ 
are independent logistic.
But this follows from the construction 
(Lemma \ref{Ljoint} and Corollary \ref{Cindep}).

\setlength{\unitlength}{0.5in}
\begin{picture}(9.5,5.9)(-0.5,-0.3)
\put(0,4){\fbox{$\phi$}}
\put(2.5,3){\fbox{$k$}}
\put(-0.6,0){\fbox{$j_1j_2 \ldots j_{h}$}}
\put(1.9,0){\fbox{$kj_1j_2 \ldots j_{h-1}$}}
\put(0.52,3.93){\line(5,-2){1.9}}
\put(0.3,3.8){\line(1,-1){0.5}}
\put(0,3.8){\line(-1,-1){0.5}}
\put(0.15,3.7){\line(0,-1){0.2}}
\put(0.15,3.3){\line(0,-1){0.2}}
\put(0.15,2.9){\line(0,-1){0.2}}
\put(0.15,2.5){\line(0,-1){0.2}}
\put(0.15,2.1){\line(0,-1){0.2}}
\put(0.15,1.7){\line(0,-1){0.2}}
\put(0.15,1.3){\line(0,-1){0.2}}
\put(0.15,0.9){\line(0,-1){0.2}}
\put(2.8,2.8){\line(1,-1){0.5}}
\put(2.5,2.8){\line(-1,-1){0.5}}
\put(2.65,2.7){\line(0,-1){0.2}}
\put(2.65,2.3){\line(0,-1){0.2}}
\put(2.65,1.9){\line(0,-1){0.2}}
\put(2.65,1.5){\line(0,-1){0.2}}
\put(2.65,1.1){\line(0,-1){0.2}}
\put(2.65,0.7){\line(0,-1){0.2}}
\put(5,4){\fbox{$l$}}
\put(7.5,3){\fbox{$\phi$}}
\put(4.4,0){\fbox{$l_1j_1 \ldots j_{h}$}}
\put(6.9,0){\fbox{$j_1j_2 \ldots j_{h-1}$}}
\put(5.52,3.93){\line(5,-2){1.9}}
\put(5.3,3.8){\line(1,-1){0.5}}
\put(5,3.8){\line(-1,-1){0.5}}
\put(5.15,3.7){\line(0,-1){0.2}}
\put(5.15,3.3){\line(0,-1){0.2}}
\put(5.15,2.9){\line(0,-1){0.2}}
\put(5.15,2.5){\line(0,-1){0.2}}
\put(5.15,2.1){\line(0,-1){0.2}}
\put(5.15,1.7){\line(0,-1){0.2}}
\put(5.15,1.3){\line(0,-1){0.2}}
\put(5.15,0.9){\line(0,-1){0.2}}
\put(7.8,2.8){\line(1,-1){0.5}}
\put(7.5,2.8){\line(-1,-1){0.5}}
\put(7.65,2.7){\line(0,-1){0.2}}
\put(7.65,2.3){\line(0,-1){0.2}}
\put(7.65,1.9){\line(0,-1){0.2}}
\put(7.65,1.5){\line(0,-1){0.2}}
\put(7.65,1.1){\line(0,-1){0.2}}
\put(7.65,0.7){\line(0,-1){0.2}}
\put(6.1,-0.5){$\dE_{h,l}$}
\put(1.1,-0.5){$\dE_h$}
\end{picture}

\vspace{0.3in}

\centerline{
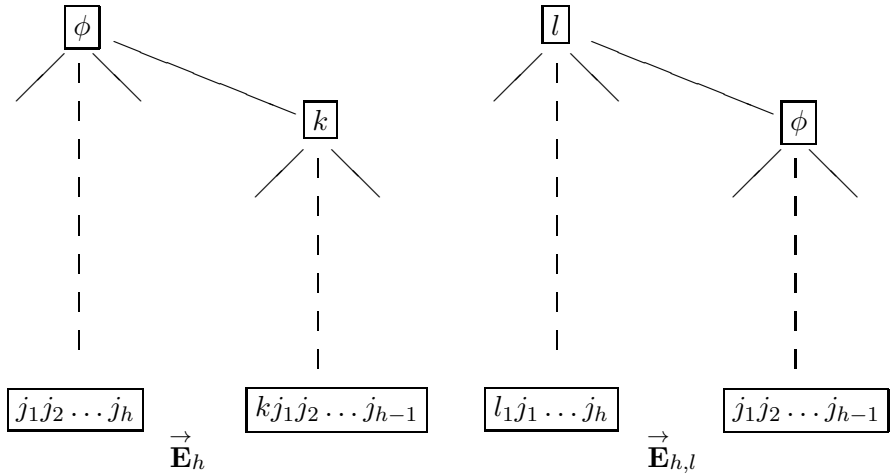
{\bf Figure 5.}
}

\subsection{The bi-infinite tree}
\label{sec-CI}
Recall $\lambda$ denotes the distribution of the edge-weights
$(W(e), e \in \bE)$ of the PWIT.
So $\lambda$ is a probability measure on $\bW = (0,\infty)^{\bE}$.
To incorporate a distinguished neighbor of $\phi$ we extend to the
state space $\bW \times \{1,2,3,\ldots \}$
and introduce the $\sigma$-finite measure
$\lambda \times {\tt count}$.
We now describe an equivalent way of representing this structure,
which turns out to be convenient for certain calculations
with $\Mopt$.

Take two copies, $\bT^+$ and $\bT^-$ say, of the PWIT,
and write their vertices as $+v$ and $-v$.
Then construct a new ``bi-infinite" tree $\bT^{\lra}$ by joining the roots $+\phi$ and
$- \phi$ of $\bT^+$ and $\bT^-$ 
via a distinguished edge $(+\phi, -\phi)$.
Write $\bE^{\lra}$ for its edge-set.
Let edge-weights on the edges of each of $\bT^+$ and $\bT^-$ be distributed as
in the PWIT, independently for the two sides of $\bT^{\lra}$.
Then define a $\sigma$-finite measure $\lambda^{\lra}$ on
$\bW^{\lra}:= (0,\infty)^{\bE^{\lra}}$
by specifying that
the weight $W(-\phi,+\phi)$ on the distinguished edge should have
``distribution" uniform on $(0,\infty)$, independent of the other edge-weights.

{\footnotesize 
\setlength{\unitlength}{0.4in}
\begin{picture}(12,9.8)(-6,-5)
\put(-1.2,-0.1){\fbox{$-\phi$}}
\put(-0.4,0){\line(1,0){1.1}}
\put(0.8,-0.1){\fbox{$+\phi$}}
\put(1.3,-0.3){\line(1,-1){1.4}}
\put(1.5,0){\line(1,0){1.2}}
\put(1.3,0.3){\line(1,1){1.4}}
\put(3.47,-2.3){\line(1,-1){1.2}}
\put(3.47,-2.0){\line(1,0){1.2}}
\put(3.47,-2.15){\line(2,-1){1.2}}
\put(3.47,-0.3){\line(2,-1){1.2}}
\put(3.47,0){\line(1,0){1.2}}
\put(3.47,0.3){\line(2,1){1.2}}
\put(3.47,2.0){\line(1,0){1.2}}
\put(3.47,2.15){\line(2,1){1.2}}
\put(3.47,2.3){\line(1,1){1.2}}
\put(2.8,-2.1){\fbox{+1}}
\put(2.8,-0.1){\fbox{+2}}
\put(2.8,1.95){\fbox{+3}}
\put(4.7,-4.0){\fbox{+11}}
\put(4.7,-3.0){\fbox{+12}}
\put(4.7,-2.1){\fbox{+13}}
\put(-0.1,-0.3){0.8}
\put(1.8,-0.3){2.1}
\put(1.8,-1.5){0.5}
\put(1.8,0.5){3.2}
\put(3.7,-3.4){1.8}
\put(3.7,-2.7){1.9}
\put(3.7,-2.3){3.7}
\put(4.7,-1.1){\fbox{+21}}
\put(4.7,-0.1){\fbox{+22}}
\put(4.7,0.9){\fbox{+23}}
\put(3.7,-0.3){2.8}
\put(3.7,-0.9){0.9}
\put(3.7,0.28){4.1}
\put(-1.3,-0.3){\line(-1,-1){1.4}}
\put(-1.3,0){\line(-1,0){1.2}}
\put(-1.3,0.3){\line(-1,1){1.4}}
\put(-3.3,-2.3){\line(-1,-1){1.2}}
\put(-3.3,-2.0){\line(-1,0){1.2}}
\put(-3.3,-2.15){\line(-2,-1){1.2}}
\put(-3.3,-0.3){\line(-2,-1){1.2}}
\put(-3.3,0){\line(-1,0){1.2}}
\put(-3.3,0.3){\line(-2,1){1.2}}
\put(-3.3,2.0){\line(-1,0){1.2}}
\put(-3.3,2.15){\line(-2,1){1.2}}
\put(-3.3,2.3){\line(-1,1){1.2}}
\put(-3.2,-2.1){\fbox{$-1$}}
\put(-3.2,-0.1){\fbox{$-2$}}
\put(-3.2,1.95){\fbox{$-3$}}
\put(-5.3,-4.0){\fbox{$-11$}}
\put(-5.3,-3.0){\fbox{$-12$}}
\put(-5.3,-2.1){\fbox{$-13$}}
\put(-2.2,-0.3){1.2}
\put(-2.2,-1.4){0.9}
\put(-2.2,0.5){2.5}
\put(-4.3,-3.5){2.2}
\put(-4.3,-2.7){2.6}
\put(-4.3,-2.3){2.7}
\put(-5.3,-1.1){\fbox{$-21$}}
\put(-5.4,-1.1){\line(-1,-1){0.7}}
\put(-5.4,-1.0){\line(-1,0){0.7}}
\put(-5.4,-0.9){\line(-1,1){0.7}}
\put(-5.3,-0.1){\fbox{$-22$}}
\put(-5.3,0.9){\fbox{$-23$}}
\put(-4.3,-0.3){1.7}
\put(-4.3,-1.0){1.1}
\put(-4.3,0.3){3.3}
\end{picture}
}

{\bf Figure 6.}  
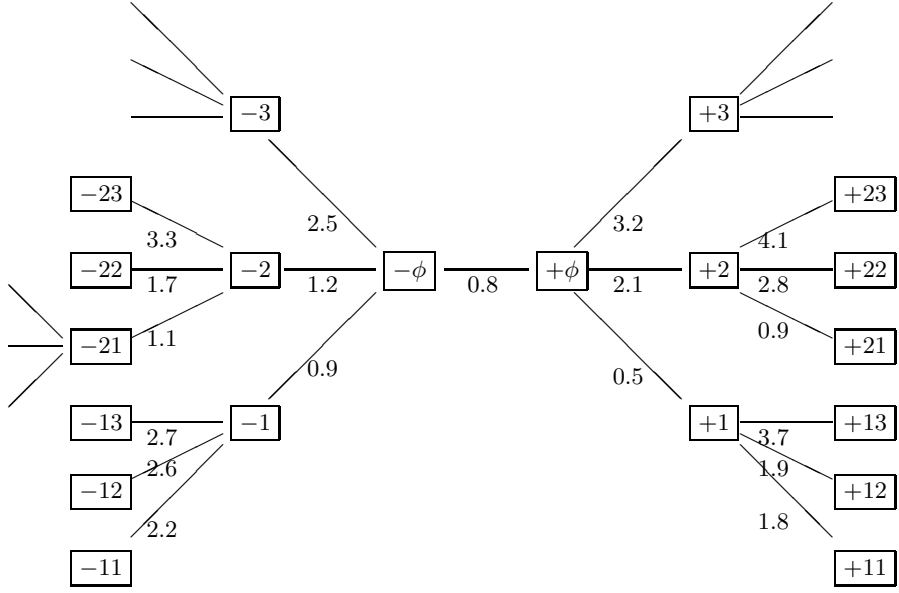
Part of a realization $\bw$ of 
edge-weights on the bi-infinite tree $\bT^{\lra}$.

\vspace{0.2in}
\noindent
The point of this construction is that there is a natural bijection between
$\bW \times \{1,2,3,\ldots\}$ and $\bW^{\lra}$
which takes $\lambda \times {\tt count}$ to $\lambda^{\lra}$.
In fact, if $k$ is the distinguished vertex of $\bT$, then
relabel vertices according to the rules

relabel $k$ as $- \phi$

relabel $j$ as $+j$ for $j \leq k-1$ and as $+(j-1)$ for $j \geq k+1$

relabel descendants accordingly.

\vspace{0.07in}
\noindent
This relabeling induces a map $\psi: 
\bW \times \{1,2,3,\ldots\} \to \bW^{\lra}$
which is invertible.
See figure 6.
Checking that $\psi$ maps $\lambda \times {\tt count}$ to $\lambda^{\lra}$ reduces
to the following easy lemma.
\begin{Lemma}
Write $\Delta := \{(x_i): 0<x_1<x_2< \ldots, x_i \to \infty\}$.
Write {\tt Pois} for the probability measure on $\Delta$ which
is the distribution of the Poisson process of rate $1$.
Consider the map
$\chi: \Delta \times \{1,2,3,\ldots\} \to \Delta \times (0,\infty)$
which takes $((x_i),k)$ to $((x_i, i \neq k), x_k)$.
Then
$\chi$ maps {\tt Pois} $\times$ {\tt count}
to {\tt Pois} $\times$ {\tt Leb},
where {\tt Leb} is Lebesgue measure on $(0,\infty)$.
\end{Lemma}
{\em Remark.}
There is a natural map ``reflect" from $\bW^{\lra}$ to $\bW^{\lra}$, induced
by the bijection of vertices $+v \lra -v$.
The map $\theta: \bW \times \{1,2,3,\ldots\} \to \bW \times \{1,2,3,\ldots,\}$ featuring in the definition of spatial invariance
(where we now ignore matchings)
is related to $\psi$ as indicated in the diagram.

\begin{center}
$\begin{array}{ccccl}
\bW \times \{1,2,\ldots\} &&
\stackrel{\psi}{\to}&&\bW^{\lra}\\
&&&&\\
\downarrow \theta &&&&\downarrow \mbox{reflect}\\
&&&&\\
\bW \times \{1,2,\ldots\} &&
\stackrel{\psi^{-1}}{\leftarrow}&&\bW^{\lra}
\end{array}$
\end{center}

\noindent
In \cite{me60},
our Theorem \ref{Told}
was stated in terms of matchings on $\bT^{\lra}$ instead of $\bT$.
In this paper we switched to using $\bT$ as the basic limit
structure for two reasons:\\
(i) on $\bT$ we can define a random matching $\MM$ using probability
distributions (instead of having to introduce $\sigma$-finite 
measures);\\
(ii) the definition of spatial invariance
(which in either setting involves $\sigma$-finite measures)
is simpler on $\bT$ than on $\bT^{\lra}$.\\
In section \ref{sec-rec} we reconcile the definitions.

The relabeling used to define $\psi$ can also be used to define
a family
$(X(\ef), \ef \mbox{ a directed edge of } \bT^{\lra})$
in terms of the $X(\ef)$ on the directed edges of the PWIT,
constructed in Lemma \ref{Ljoint}.
It is easy to check that the joint distribution of
$(W(e), X(\ef); e, \ef \mbox{ edges of } \bT^{\lra})$
thus obtained is the same as if we applied the construction
in Lemma \ref{Ljoint} to
$(W(e), e \in \bE^{\lra})$, replacing 
$\dE_h$ in the construction by
\begin{equation}
\stackrel{\lra}{bE}_h
 = \{\ef = (+\bv, +\bv j): \ |\bv| = h, j \geq 1 \}
\cup   \{\ef = (-\bv, -\bv j): \ |\bv| = h, j \geq 1 \}
. \label{Elra}
\end{equation}
Then the matching $\Mopt$ can be defined on $\bT^{\lra}$ in the same way
as on $\bT$:
\[ 
e \mbox{ is an edge of $\Mopt$ iff }
W(e) < X(\ef) + X(\eb)
. \]
Lemma \ref{Llam} below shows how working on the bi-infinite
tree is useful for calculations.
Informally, Lemma \ref{Llam} describes the distribution
of $\Mopt$ as seen from a typical edge 
in $\Mopt$, and exhibits a conditional independence property
for the restrictions of $\Mopt$ to the two sides of the tree
determined by that edge.

On the PWIT define
$ X^{\downarrow} = \min_{i \geq 1} (W(\phi,i) - X(\phi,i)) $,
and write $\nu_x$ for the conditional distribution of
the family
$(W(e), e \in \bE; X(\ef), \ef \in \dE, \ef
\mbox{ directed away from } \phi)$
given $X^{\downarrow} = x$.
Returning to the bi-infinite tree,
write $\lambda^1$ for the measure obtained by restricting
$\lambda^{\lra}$ to the set
$\{W(-\phi,+\phi) < X(-\phi,+\phi) + X(+\phi,-\phi)\}$.
So under $\lambda^1$, in $\Mopt$ the vertex $+\phi$ is a.s. matched
with vertex $- \phi$.
\begin{Lemma}
\label{Llam}
$\lambda^1$ is a probability measure.
Under $\lambda^1$ we have:\\
(i) the joint density of $(W(-\phi,+\phi),X(-\phi,+\phi),X(+\phi,-\phi))$
at $(w,x_1,x_2)$ is
equal to
$f(x_1)f(x_2)1_{(0<w<x_1+x_2)}$,
where $f$ is the logistic density; \\ 
(ii) conditional on
$(W(-\phi,+\phi),X(-\phi,+\phi),X(+\phi,-\phi)) = (w,x_1,x_2)$,
the distribution of
the family
\[ (W(e), e \in \bE^+; X(\ef), \ef \in \dE^+, \ef
\mbox{ directed away from } + \phi) \]
is the image of $\nu_{x_1}$ under the natural embedding
$\bT \to \bT^+ \subset \bT^{\lra}$;
the distribution of
the family
\[ (W(e), e \in \bE^-; X(\ef), \ef \in \dE^-, \ef
\mbox{ directed away from } - \phi) \]
is the image of $\nu_{x_2}$ under the natural embedding
$\bT \to \bT^- \subset \bT^{\lra}$;
and these two families are conditionally independent.
\end{Lemma}
Note that because $+\phi$ is matched to $-\phi$, we have
$X(+j,+\phi) = W(+\phi,-\phi) - X(+\phi,-\phi)$
and then one can recursively construct $X(\ef)$ for $\ef$
directed {\em toward} $(-\phi,+\phi)$.
So the prescription in Lemma \ref{Llam} is enough to
specify the joint distribution of all the $X(\ef)$ and hence
of $\Mopt$, under $\lambda^1$.

{\em Proof of Lemma \ref{Llam}.}
The joint density has the form stated in (i) by construction;
and so its total mass equals
$\int \int (x_1 + x_2)^+ \ f(x_1)f(x_2) \ dx_1 dx_2$
which equals $1$ by Lemma \ref{L1}.
Next, by the construction based on 
$(X(\ef), \ef \in 
\stackrel{\lra}{bE}_h
)$ at (\ref{Elra}), we see that under $\lambda^{\lra}$
the families
\[ \{W(e), e \in \bE^+; X(\ef), \ef \in \dE^+, \ef
\mbox{ directed away from } + \phi\} 
\ \cup \\{ X(-\phi,+\phi)\\} \]
and
\[ \{W(e), e \in \bE^-; X(\ef), \ef \in \dE^-, \ef
\mbox{ directed away from } - \phi\}
\ \cup \\{ X(+\phi,-\phi)\\} \]
are independent of each other and of $W(-\phi,+\phi)$.
Since $\lambda^1$ is defined by an event depending only
on 
$\{  X(+\phi,-\phi),X(-\phi,+\phi),W(+\phi,-\phi) \}$,
we obtain the desired conditional independence property.
Each family under $\lambda^{\lra}$ is the distributed as
the image of the corresponding family on the PWIT
(making $X(-\phi,+\phi)$ correspond to $X^\downarrow$),
and so by independence under $\lambda^{\lra}$ the conditional
distribution under $\lambda^1$ depends only on $x_1$
(resp. $x_2$).
$\Box$

Recall the map
$\psi: 
\bW \times \{1,2,3,\ldots\} \to \bW^{\lra}$
which takes $\lambda \times {\tt count}$ to $\lambda^{\lra}$.
The inverse image of the event
$\{W(-\phi,+\phi) < X(-\phi,+\phi) + X(+\phi,-\phi)\}$
is
\[ \psi^{-1}
\{W(-\phi,+\phi) < X(-\phi,+\phi) + X(+\phi,-\phi)\}
\ = \ 
\{\Mopt(\phi) \mbox{ is distinguished} \} . \]
Thus the inverse image of the probability measure $\lambda^1$
is $\lambda \times {\tt count}$ restricted to 
$\{\Mopt(\phi) \mbox{ is distinguished} \} $. 
Then after un-distinguishing this vertex, we are left with
probability distribution $\lambda$ on the PWIT.

So in summary,
we have established the following ``relabeling principle".
\begin{quote}
Given $(Z(\ef))$ and hence $\Mopt$ on the PWIT, map 
the whole structure
to the bi-infinite tree by relabeling $(\phi,\bMopt(\phi)$ as $(+\phi,-\phi$
and relabeling other vertices accordingly; then the
resulting
distribution on the bi-infinite tree is $\lambda^1$.
\end{quote}
To see why this is useful, let us give a quick second proof of
Proposition \ref{P8}(a).
The distribution of $W(\phi,\bMopt(\phi))$ in the PWIT
is the same,
by the relabeling principle,
as the distribution on the bi-infinite tree of
$W(-\phi,+\phi)$ under $\lambda^1$.
Then by Lemma \ref{Llam}(i)
\begin{eqnarray*}
 \lambda^1\{ W(-\phi,+\phi) \in [w,w+dw] \} /dw
&=& \int \int_{x_1+x_2 > w} f(x_1)f(x_2) \ dx_1dx_2
\\ 
&=& P(X_1+X_2>w) 
\end{eqnarray*}
for independent logistic $X_1, X_2$.

\subsection{Calculating with the bi-infinite tree}
\label{sec-CI2}
On the bi-infinite tree define
\[ C^* := \{
+ \phi = \arg \min_{y \sim +I} \, \! ^{[2]} (W(+I,y) - X(+I,y))
\} \]
where
\[ I = \arg \min_{i \geq 1} (W(+\phi,+i) - X(+\phi,+i)) . \]
\begin{eqnarray*}
\FF^{+} &=& \sigma(X(\ef),W(e): e, \ef \mbox{ edges of } \bT^{+})\\
\FF^{-} &=& \sigma(X(\ef),W(e): e, \ef \mbox{ edges of } \bT^{-})\\
\FF^{\phi} &=& \sigma (X(+\phi,-\phi),X(-\phi,+\phi),W(+\phi,-\phi)) .
\end{eqnarray*}
\begin{Lemma}
\label{lamD}
$\lambda^1 \{C^{*c}|\FF^{-},\FF^{\phi}\} = 
 g(X(-\phi,+\phi),W(+\phi,-\phi)-X(+\phi,-\phi)) $ 
for $g$ defined in Lemma \ref{LXI}.
\end{Lemma}
{\em Proof.}
$C^*$ is $\FF^+$-measurable, so by the conditional independence
assertion of Lemma \ref{Llam}
we have
\[ \lambda^1 \{C^{*c}|\FF^{-},\FF^{\phi}\} = 
 \lambda^1 \{C^{*c}|\FF^{\phi}\} . \] 
Thus we have to show that
\begin{equation}
 \lambda^1 \{C^{*c} |
(W(-\phi,+\phi),X(-\phi,+\phi),X(+\phi,-\phi)) = (w,x_1,x_2)
\} = g(x_1,w-x_2) . \label{lCW}
\end{equation}
Under this conditioning, Lemma \ref{Llam} implies that
the family $(W(e), e \in \bE^+; X(\ef), \ef \in \dE^+, \ef
\mbox{ directed away from } + \phi)$
is distributed as the image of 
the family $(W(e), e \in \bE; X(\ef), \ef \in \dE, \ef
\mbox{ directed away from } \phi)$
conditioned on $\{X^\downarrow = x_1\}$.
By definition of $g(a,b)$ in Lemma \ref{LXI},
\begin{eqnarray*}
g(x_1,w-x_2) \! &=&
\lambda^1 \{ W(+\phi,+I) - (w-x_2) >
  \min_{y \sim +I} \, \! ^{[2]} (W(+I,y) - X(+I,y))
\\
&&  |
(W(-\phi,+\phi),X(-\phi,+\phi),X(+\phi,-\phi)) = (w,x_1,x_2)
\} .
\end{eqnarray*}
But under this conditioning
\begin{eqnarray*}
W(+\phi,+I) - (w-x_2) &=&
W(+\phi,+I) - (W(+\phi,-\phi) - X(+\phi,-\phi))\\
&=& W(+I,+\phi) - X(+I,+\phi)
\end{eqnarray*}
by (\ref{Xrec}), because under $\lambda^1$ 
the vertex $+\phi$ is always matched to $-\phi$.
Substituting into (\ref{lCW}), we see that the event in (\ref{lCW})
is precisely the event $C^{*c}$, as required.

{\em Proof of Lemma \ref{LBB*}.}
The relabeling principle shows that $P(B^*|B)$ can be rewritten
as $\lambda^1\{C^* |C\}$, for a certain event $C$ which is
$\FF^\phi$-measurable
and such that $P(B) = \lambda^1 \{C\}$.
Now
\begin{eqnarray*}
\lambda^1\{C^{*c}|C\} \! &=& \!
E_{\lambda^1} \ 1_C 
 g(X(-\phi,+\phi),W(+\phi,-\phi)-X(+\phi,-\phi)) \mbox{ by Lemma \ref{lamD}}\\ 
&>& 0 \mbox{ if } \lambda^1\{C\} > 0 \mbox{ by Lemma \ref{LXI}}
\end{eqnarray*}
establishing Lemma \ref{LBB*}.

\subsection{Remarks on quantifying Theorem 4}
\label{sec-quant}
Recall equation (\ref{DA}):
the difference in cost between $\Mopt$ and another spatially invariant
random matching $\MM$ equals $ED1_A$.
To quantify Theorem \ref{T3}, we would like to show
\begin{quote}
if $P(\bMM(\phi) \neq \bMopt(\phi)) = \delta$
then $ED1_A \geq \eps(\delta)$
\end{quote}
with some explicit lower bound on $\eps(\delta)$:
our current ``proof by contradiction" of Theorem \ref{T3}
shows only that $\eps(\delta) > 0$.
Now it is not difficult to improve Lemma \ref{LXI} to a stronger
result giving a lower bound on 
a corresponding quantity, under the same conditioning.
Given some spatially invariant $\MM$, one can condition on the
one-sided infinite path 
$+\phi,-\phi,-J,\ldots$
in $\bT^{\lra}$
whose edges alternate between $\Mopt$ and $\MM$,
and seek to apply the improved lemma to the conditional distribution
of the edge $(+\phi,\bMM(+\phi))$.
But the difficulty with this scheme is that the definition
of $\MM$ might depend on the whole tree; we cannot argue 
that {\em a priori} the behavior of $\MM$ on $\bT^+$ and on
$\bT^-$ has some conditional independence property,
as exploited in the proof of Lemma \ref{lamD}.

\subsection{Reconciliation}
\label{sec-rec}
Here we reconcile the way definitions and results were stated
in \cite{me60} with the way they are stated in this paper.

Write $\bW^* = \bW \times \{1,2,3,\ldots\}$.
The map $\theta$ in the definition of spatial invariance
is a map $\theta: \bW^* \to \bW^*$
which preserves $\lambda \times {\tt count}$.
Given a spatially invariant random matching $\MM$, define
$\gamma: \bW^* \to [0,1]$ by
\begin{equation}
 \gamma(\bw,i) = P(\bMM(\phi) = i|W(e) = w(e) \ \forall e \in \bE) . \label{gwM}
\end{equation}
Then $\gamma$ must have the following two properties;
the first because $\MM$ is a matching, the second because of
spatial invariance.

\noindent
(i) $\sum_i \gamma(\bw,i) = 1 $.\\
(ii) $\gamma(\theta(\bw,i)) = \gamma(\bw,i)$.

\noindent
And we can write the associated cost as\\
(iii) $ EW(\phi,\bMM(\phi)) = \int_{\bW} \sum_i \ w(\phi,i) \gamma(\bw,i)
\ \lambda(d\bw) $. 

It turns out that we can reverse the argument: given a function
$\gamma(\cdot)$ satisfying (i) and (ii), one can define a
spatially invariant random matching satisfying
(\ref{gwM}), whose cost therefore satisfies (iii).
Now use the bijection $\psi$ to define
\[ g(\bw^{\lra}) = \gamma(\psi^{-1}(\bw^{\lra})) . \]
Then $g(\cdot)$ satisfies certain consistency conditions 
corresponding to (i,ii), and the cost (iii) becomes \\
 (iii*)
$\quad \int_{\bW^{\lra}} g(\bw^{\lra}) \ \lambda^{\lra}(d\bw^{\lra}) $.  \\
Thus the quantity 
$ \inf 
\{ E W(\phi, \bMM(\phi)) : \quad \MM \mbox{ spatially invariant} \} $ 
in our Theorem \ref{Told} equals the infimum of (iii*)
over functions $g$ satisfying
the consistency conditions, and that was the formulation of the
limit constant $\lim_n EA_n$ in \cite{me60} Theorem 1.

Our Theorem \ref{Told-easy} is \cite{me60} Proposition 3(a),
stated there for the optimal $n$-matchings but extending unchanged
to arbitrary spatially invariant matchings.
Our Theorem \ref{Told-hard} is obtained by combining
Proposition 3(b) and Proposition 2 of \cite{me60};
local convergence is the method of proof of those results.

\section{Variants of the random assignment problem}
\label{sec-VAR}
From our asymptotic viewpoint, the non-bipartite matching problem
-- with $n \times n$ matrix $(c(i,j))$ with $c(i,j) \equiv c(j,i)$
and $n$ even -- is the same (up to normalization convention) as the bipartite problem, because
the limit random structure and problem is exactly the same matching problem on the PWIT.
Below we describe variants where the limit random structure
or problem is different.

\subsection{Power-law densities}
\label{sec-power}
Consider the random assignment problem
\[ A_n = \min_\pi \sum_{i=1}^n c(i,\pi(i))  \]
where now the $(c(i,j))$ are i.i.d. with density $f_c(\cdot)$
satisfying
\begin{equation}
 f_c(x) \sim x^r \mbox{ as } x \downarrow 0 \label{fc}
\end{equation}
for some $0<r<\infty$.
This setting is motivated by trying to ``fit" mean-field models
to the Euclidean matching problem on random points in
$[0,1]^d$, for which the distribution of inter-point distances
satisfies (\ref{fc}) with $r=d-1$:
see \cite{HMM98} for further discussion.
In this setting, $EA_n$ will grow as order $n^{\frac{r}{r+1}}$,
and so we rescale the problem to study 
$A_n^\prime = n^{-\frac{r}{r+1}} A_n$; in other words 
\[ A^\prime_n = \min_\pi \sfrac{1}{n} \sum_{i=1}^n c^\prime(i,\pi(i))  \]
where
$c^\prime(i,j) = n^{1/(r+1)} c(i,j)$.
M\'{e}zard - Parisi (\cite{MP85} 
eq. (22,23))
use the replica method to argue
$EA^\prime_n \to \gamma_r$
where the limit constant is characterized by
\begin{eqnarray}
\gamma_r &=& (r+1) \int_{-\infty}^\infty G(l) e^{-G(l)} \ dl \label{r1}\\
G(l) &=& \sfrac{2}{r!} \int_{-l}^\infty (l+y)^r e^{-G(y)} \ dy,
\quad -\infty < l < \infty . \label{r2}
\end{eqnarray}
Let us indicate how our approach gives essentially the same result.
Underlying the connection between the matrix with i.i.d. exponential
 (mean $n$)
entries and the PWIT 
is the fact that the order statistics
$(\xi^{(n)}_1,\xi^{(n)}_2,\ldots)$
of $n$ exponential (mean $n$) r.v.'s satisfy
\begin{equation}
(\xi^{(n)}_1,\xi^{(n)}_2,\ldots)
\cd (\xi_1,\xi_2,\ldots) \label{conv2Pois}
\end{equation}
where the limit is a Poisson (rate $1$) process.
In the current setting, the order statistics $(\xi^{(n)}_i)$
of $n$ r.v.'s distributed as $c^\prime(i,j)$ satisfy (\ref{conv2Pois})
where the limit $(\xi_i)$ is now the 
inhomogeneous Poisson process of rate 
\begin{equation}
P(\mbox{ some point of $(\xi_i)$ in } [x,x+dx]) = x^r \ dx.
\label{Pinhom}
\end{equation}
Without checking details, it seems clear that 
following the method of \cite{me60} and this paper
leads to the formulas
which parallel (\ref{over1},\ref{over2})
\begin{eqnarray}
\gamma_r &=& \int_0^\infty x \cdot x^r P(X_1+X_2>x) \ dx \label{me1}\\
X & \ed & \min_i (\xi_i - X_i) \label{me2}
\end{eqnarray}
where $(\xi_i)$ is now the inhomogeneous Poisson process (\ref{Pinhom}).
The calculations below will check that (\ref{r1},\ref{r2})
and (\ref{me1},\ref{me2}) are essentially equivalent.
Writing $F(x) = P(X>x)$, (\ref{me2}) says
\begin{eqnarray}
F(x)&=&
P(\mbox{ no point of $(\xi_i,X_i)$ lies in } 
\{(z,b): \ z-b \leq x\} ) \nonumber\\
&=& \exp \left( - \int_{0}^\infty z^r F(z-x) \ dz \right) .
\label{zyx}
\end{eqnarray}
Setting $y = z-x$ gives
\[ F(x) = \exp \left( - \int_{-x}^\infty (y+x)^r F(y) \ dy \right) . \]
Setting $G(x) = - \log F(x)$ gives
\[ G(x) = \int_{-x}^\infty (y+x)^r e^{-G(y)} \ dy \]
which agrees with (\ref{r2}) up to the constant factor
$2/r!$.  
[why the discrepancy?  The $2$ is the normalization convention
mentioned in section \ref{sec-calc}.  I guess the $r!$ is
a typo at \cite{MP85} eq. (5), where the $l^r$ should be
$l^r/r!$ as it is in \cite{MP86} eq. (2.1).]
Next, (\ref{me1}) says
\begin{eqnarray*}
\gamma_r &=&
\int_0^\infty x^{r+1} dx \ 
\int_{-\infty}^\infty P(X_1 \in dy) P(X_2 \geq x-y) \\
&=& \int_{-\infty}^\infty P(X \in dy) \ 
\int_0^\infty x^{r+1} F(x-y) \ dx \\
&=& - \int_{-\infty}^\infty F(y) dy \ 
\int_0^\infty x^{r+1} F^\prime(x-y) dx \\
&& \mbox{(integrating by parts over $y$)}\\
&=&  \int_{-\infty}^\infty F(y) dy \ 
(r+1) \int_0^\infty x^r F(x-y) dx \\
&& \mbox{(integrating by parts over $x$)}\\
&=& (r+1) \int_{-\infty}^\infty F(y) (- \log F(y)) \ dy \ \ \mbox{ by (\ref{zyx})} \\
&=& 
(r+1) \int_{- \infty}^\infty G(y) e^{-G(y)} \ dy
\ \ \mbox{ for } G(y) = - \log F(y)
\end{eqnarray*}
and this is exactly (\ref{r1}).

\subsection{The combinatorial TSP}
\label{sec-TSP}
In the combinatorial (i.e. symmetric mean-field) 
traveling salesman problem (TSP)
we take a $n \times n$ matrix of ``distances" $c(i,j)$
such that
$(c(i,j), 1 \leq i < j \leq n)$ are i.i.d. exponential (mean $n$)
and $c(i,i) \equiv 0$ and $c(j,i) \equiv c(i,j)$.
We study
\[ L_n:= \min_\pi \sfrac{1}{n} \sum_i c(i,\pi(i)) \]
where $\pi$ is a cyclic permutation of $\{1,2,\ldots,n\}$.
One can mimic the heuristic argument of section \ref{sec-heur}.
That is, for the subtree $\bT^v$ rooted at $v$, a {\em tour}
is a set of doubly-infinite paths which pass once through each vertex of
$\bT^v$; an {\em almost-tour} is the variation in which the root $v$ is 
the start of one path in the set.
Write (heuristically)
\begin{eqnarray*}
(*) \quad \quad \quad X_{v} &=&
\mbox{ cost of optimal tour on } \bT^v\\
 &-& \mbox{ cost of optimal almost-tour on } \bT^v.
\end{eqnarray*}
One can continue the parallel heuristic
to argue that the solution format 
of the combinatorial TSP is very similar to that for the
random assignment problem given by (\ref{over1},\ref{over2}):
\begin{equation}
\lim_n EL_n = \int_0^\infty z P(X_1+X_2 \geq z) \ dz 
\label{TSP2}
\end{equation}
where the distribution of $X, X_i$ is determined by
\begin{equation}
 X \ed \min_i \, \! ^{[2]} (\xi_i - X_i)
\mbox{ where $(\xi_i)$ is a Poisson($1$) process;}
\label{TSP1}
\end{equation}
where $\min^{[2]}$ denotes the second-smallest value
(so the only difference between the two solutions is the
replacement of $\min$ by $\min^{[2]}$).
By copying the proof of Lemma \ref{Llogistic}
one can progress toward an explicit solution.
Writing $\bar{F}(y) = P(X > y)$,
(\ref{TSP1}) says (in the notation of Lemma \ref{Llogistic})
\begin{eqnarray*}
\bar{F}(y)&=&
 P(\mbox{ $0$ or $1$  points of $\PP$ in $\{(z,x):\ z-x \leq y\}$})\\
&=& \left( 1+\int_{-y}^{\infty} \bar{F}(u) du\right)
\exp \left(- \int_{-y}^{\infty} \bar{F}(u) du \right) .
\end{eqnarray*}
Now write
$G(y) = \int_{-y}^{\infty} \bar{F}(u) du$
so that
the fixed-point equation (\ref{TSP1}) becomes 
the fixed-point equation  
\begin{equation}
 G^\prime(y) = (1+G(-y)) \exp(-G(-y)) . \label{TSP3}
\end{equation}
To rewrite (\ref{TSP2}) in terms of $G$,
set $z = x_1+x_2$.
\begin{eqnarray*}
 \int_0^\infty z P(X_1+X_2 > z) \ dz 
&=&
\int_0^\infty z \ dz \ \int_{-\infty}^\infty
F^\prime(x_2) \bar{F}(z-x_2) \ dx_2\\
&=&
\int_0^\infty z \ dz \ \int_{-\infty}^\infty
F^\prime(z-x_1) \bar{F}(x_1) \ dx_1\\
&=& \int_{-\infty}^\infty \bar{F}(x_1) \ dx_1 \ \int_0^\infty
z F^\prime(z-x_1) \ dz\\
&=& \int_{-\infty}^\infty \bar{F}(x_1) \ dx_1 \ \int_0^\infty
F(z-x_1) \ dz\\
&=& \int \int_{x_1+x_2>0} \bar{F}(x_1)F(x_2) \ dx_1dx_2\\
&=& \int_{-\infty}^\infty  \ - G^\prime(x_1)G(-x_1) \ dx_1 .
\end{eqnarray*}
Now use (\ref{TSP3}) to see that (\ref{TSP2}) becomes
\begin{equation}
\lim_n EL_n = \int_{-\infty}^\infty G(x)(1+G(x))e^{-G(x)} \ dx .
\label{TSP4}
\end{equation}
So we can rewrite our solution (\ref{TSP2},\ref{TSP1})
as (\ref{TSP3},\ref{TSP4}),
which is the solution given by
Krauth - M{\'e}zard \cite{KM87}
using the cavity method.
Earlier work \cite{MP86} using the replica method gave more complicated
formulas.
As in section \ref{sec-power} one can consider the more general setting
(\ref{fc}) and, after some calculus, our solution (with $\xi_i$ now 
defined via (\ref{Pinhom})) again coincides
with the general-$r$ expression in \cite{KM87}.
Our probabilistic expressions (\ref{TSP2},\ref{TSP1};\ref{over1},\ref{over2})
perhaps makes the mathematical similarities between the heuristic solutions 
of TSP and the random assignment problems more 
visible than do the analytic expressions in the physics literature.
In contrast to the random assignment case, the solution of identity (\ref{TSP3})
seems to have no simple explicit formula.
Numerically solving (\ref{TSP3},\ref{TSP4}) gives a limit constant of about $2.04$
\cite{KM87},
in agreement with Monte Carlo simulations of the finite-$n$
combinatorial TSP,
and the agreement persists for larger values of $r$;
see Percus - Martin \cite{PM99} for a recent review.

Returning to the bipartite setting, one can define a $k$-assignment
problem: study
\[ A_n^{(k)} := \min_S \sfrac{1}{nk} \sum_{(i,j) \in S} c_{ij} \]
where $S$ denotes an edge-set such that
$|\{j: \ (i,j) \in S\}| = k \ \forall i$
and
$|\{i: \ (i,j) \in S\}| = k \ \forall j$.
A very similar heuristic argument indicates that
$\lim_n EA_n^{(k)}$ should be given by the analog of
(\ref{TSP2},\ref{TSP1})
with $k$'th minimum in place of second-minimum.

\subsection{Gibbs measures on assignments}
\label{sec-Gibbs}
Fix a parameter $0<\lambda<\infty$.
Take $(c(i,j), 1 \leq i,j \leq n)$
i.i.d. exponential (mean $n$).
Define a non-uniform random permutation $\Pi_n$ of $\{1,2,\ldots,n\}$ by
\begin{equation}
 P(\Pi_n = \pi|\mbox{ all the } c(i,j)) \propto 
\exp( - \lambda \sum_i c(i,\pi(i))) . \label{Gibbs-form}
\end{equation}
The statistical physics approach to the random assignment
problem (and other combinatorial optimization problems over random data)
is to first study such {\em Gibbs measures}
and then take $\lambda \to \infty$ limits.

It seems plausible that our method extends to the study of
the $n \to \infty, \ \lambda$ fixed, limit behavior of this
Gibbs measure, though it may be technically challenging to make
these ideas rigorous.
In brief, regard the right side of (\ref{Gibbs-form}) as specifying
a random matching $\Pi_\infty$ on the PWIT,
and seek the density $h_\lambda(x)$ of cost-$x$ edges in the
matching $\Pi_\infty$.
To mimic the heuristic argument of section \ref{sec-heur},
then regard the right side of (\ref{Gibbs-form}) as specifying
a ``matching or almost-matching" $\Pi^+$ on $\bT^+$,
and define
\[ X_\phi = \frac{P(\Pi^+ \mbox{ is an almost-matching })}
{P(\Pi^+ \mbox{ is a matching })} . \]
One can argue heuristically that this $X_\phi$ should
satisfy the distributional identity
\begin{equation}
X \ed 
 \left( \sum_{i=1}^\infty e^{- \lambda \xi_i} X_i \right)^{-1}  \label{Gibbs-id}
\end{equation}
for Poisson (rate $1$) $\xi_i$,
and that the desired limit density is
\begin{equation}
 h_\lambda(x) = 1 - 
E \left( \frac{1}{1+X_1X_2e^{-\lambda x}} \right) . \label{hGibbs}
\end{equation}
Similar ideas are implicit in
Talagrand \cite{tal01}, who obtains some rigorous results
for small $\lambda$.
M\'{e}zard - Parisi (\cite{MP85} 
eq. (18-20)) derive a non-rigorous solution in a different form;
presumably (cf. sections \ref{sec-power} and \ref{sec-TSP}) the two forms are equivalent
but we have been unable to verify this.
It is easy to verify that, writing $X = e^{-\lambda U}$ in (\ref{Gibbs-id}),
expressions (\ref{Gibbs-id},\ref{hGibbs}) are consistent in the
$\lambda \to \infty$ limit with (\ref{over1},\ref{over2}).
Formalizing our approach in the Gibbs setting seems more
challenging than in previous settings, for the following reason.
A central part of \cite{me60} is showing (Proposition 2) that
given an almost-complete matching,
one can construct a complete matching for small extra cost.
Proving an analogous result in the Gibbs setting will presumably
require the same type of technical estimates needed in \cite{tal01}.

\section{Final remarks}
\label{sec-FR}
\subsection{The AEU property in random optimization}
\label{sec-AEU}
One can define the AEU property in quite general
``optimization over random data" settings.
Consider for instance the Euclidean 
TSP involving $n$ i.i.d. points on $[0,1]^2$.
Writing $L(T_n)$ for the length of the shortest tour,
we have \cite{BHH59} 
$EL(T_n) \sim c n^{1/2}$
for constant $c$.
The AEU property would be:
for non-optimal tours $T_n^\prime$,
if $E n^{-1} \# \{e: e \in T_n \setminus T_n^\prime \} \geq \delta$
then
$\liminf_n n^{-1/2} (EL(T_n^\prime) - EL(T_n)) \geq \eps(\delta)$
where $\eps(\delta) > 0$.
Whether or not the AEU property holds 
for the Euclidean TSP is a challenging question.
Conceptually, the AEU property seems an interesting way of
classifying optimization problems from the average-case
viewpoint and seems worthy of further study.  Note it is an intrinsic property of the problem,
not of any particular algorithm to solve the problem.

\subsection{The exact conjecture}
Our proofs are purely asymptotic, so do not shed light on
the exact conjecture (\ref{exact}).
Note (cf. proof of Proposition \ref{P8}) that $h(x)$ represents the asymptotic chance than a cost-$x$ edge
gets into the optimal matching, so in particular a cost-$0$ edge
has asymptotic chance $h(0) = 1/2$ to be in the optimal matching, as suggested
by previous work (\cite{olin}, \cite{LW00} sec. 4.1).

It is very natural to conjecture that
$\var A_n \sim \sigma^2/n$ for some
$0<\sigma<\infty$,
and that rescaled $A_n$ has a Normal limit:
this is supported by Monte Carlo simulation \cite{HMM98}.
One can define a candidate value $\tilde{\sigma}$ in terms of
the optimal matching on the infinite tree, which 
(roughly speaking) reflects {\em local} dependence, but it is
not apparent even heuristically whether $\tilde{\sigma}$
should coincide with $\sigma$.

Parisi \cite{parisi98} observes that 
$EA_n^2 \approx EA_n + 1$ for small $n$.
The implicit conjecture is an amusing instance of how numerics
can mislead; what's really going on is 
\[ EA_n^2 \approx (EA_n)^2 \mbox{ for large } n \]
and $\pi^2/6$ just happens to be close to the solution of
$z^2 = z+1$.

\subsection{Frieze's $\zeta(3)$ result}
It is intriguing that the limit $\zeta(3)$ arises in the
somewhat analogous problem of the minimum-weight spanning tree
on the complete graph with i.i.d. edge-weights
(Frieze \cite{fri85}).
That problem can also be studied using the PWIT
(Aldous \cite{me49}).
The density function $\hat{h}(x)$ with mean $\zeta(3)$ analogous to $h(x)$ is
\[ \hat{h}(x) = 
\sfrac{1}{2} (1 - q^2(x)) \]
where $q(x) = 0$ on $0 \leq x \leq 1$ and for $x>1$ 
is the non-zero solution of
\[ 1 - q(x) = \exp(-xq(x)) . \]
Thus we see the same structure: calculations of asymptotic quantities
involve a fixed-point identity.
It is natural to speculate that some other problem has solution
$\zeta(4)$.

\vspace{0.2in}

{\em Acknowledgements.}
I thank Greg Sorkin, Boris Pittel and Michel Talagrand for
discussions and encouragement.

\end{document}